\documentclass[]{amsart} 

\usepackage{amsmath} 
\usepackage{latexsym}
\usepackage{amssymb}
\ifx\pdftexversion\undefined
  \usepackage[dvips]{graphics}
\else
  \usepackage[pdftex]{graphics}
\fi
\newtheorem{Proposition}{Proposition}
  \newtheorem{Remark}{Remark}
  \newtheorem{Corollary}[Proposition]{Corollary}
  \newtheorem{Lemma}[Proposition]{Lemma}
   \newtheorem{Theorem}[Proposition] {Theorem}
\newtheorem{Condition}[Proposition]{Condition}

\newtheorem{Definition}[Proposition]{Definition}
\newtheorem{Note}[Proposition]{Note}

\def\mb{\mathbf}
\def\z{\noindent}
\def\CC{\mathbb{C}}
\def\RR{\mathbb{R}}
\def\NN{\mathbb{N}}
\def\ZZ{\mathbb{Z}}
\def\lap{\mathcal{L}}
\def\bor{\mathcal{B}}
\def\phi{\varphi}
\def\bfk{\mathbf{k}}
\def\bflam{\boldsymbol{\lambda}}
\def\CC{\mathbb{C}}
 \def\RR{\mathbb{R}}
 \def\NN{\mathbb{N}}
\def\ZZ{\mathbb{Z}}
\def\bor{\mathcal{B}}
\def\lap{\mathcal{L}}

\def\bk{\mathbf{k}}
\def\bfC{\mathbf{C}}
\def\bj{\mathbf{j}}
\def\bfm{\mathbf{m}}

\def\bfbet{\boldsymbol{\beta}}
\def\smallmo{\ensuremath{\,\rule{0.05em}{0.7em}\rule{0.5em}{0.05em}\rule{0.05em}{0.7em}}\,}
\def\largemo{\ensuremath{\,\rule{0.05em}{0.7em}\rule[0.7em]{0.5em}{0.05em}\rule{0.05em}{0.7em}}\,}

\def\anymo{\ensuremath{\,\rule{0.05em}{0.7em}\rule{0.55em}{0.05em}\hskip-0.55em
\rule[0.7em]{0.55em}{0.05em}\hskip-0.05em\rule{0.05em}{0.7em}}\,}

\def\smalltr{\,\ensuremath{\rule{0.05em}{0.7em}\rule{0.5em}{0.05em}\rule{0.05em}{0.7em}\rule{0.5em}{0.05em}\rule{0.05em}{0.7em}}\,}

\def\largetr{\,\ensuremath{\rule{0.05em}{0.7em}\rule[0.7em]{0.5em}{0.05em}\rule{0.05em}{0.7em}s
\rule[0.7em]{0.5em}{0.05em}\rule{0.05em}{0.7em}}\,}

\def\anytr{\,\ensuremath{\rule{0.05em}{0.7em}\rule{0.55em}{0.05em}\hskip-0.55em
\rule[0.7em]{0.55em}{0.05em}\hskip-0.05em\rule{0.05em}{0.7em}
\rule{0.55em}{0.05em}\hskip-0.55em
\rule[0.7em]{0.55em}{0.05em}\hskip-0.05em\rule{0.05em}{0.7em}}\,}
\def\spacetr{\stackrel{\approx}{\mathcal{T}}}
\def\erm{\mathrm{e}}
\def\mag{\,\mathrm{mag}}
\def\dom{\,\mathrm{dom}}

\def\drm{\mathrm{d}}

\begin{document}

\begin{abstract}
  
  Transseries in the sense of \'Ecalle are constructed using a
  topological approach. A general contractive mapping principle is
  formulated and proved, showing the closure of transseries under a
  wide class of operations.  
  
  In the second part we give an overview of results and methods
  reconstruction of actual functions and solutions of equations from
  transseries by generalized Borel summation with in ordinary and
  partial differential and difference equations.

\end{abstract}

\author{O.  Costin } 
\address{Department of Mathematics\\ Rutgers University
\\ Piscataway, NJ 08854-8019}

\title{Topological construction of transseries and introduction to
  generalized Borel summability}\gdef\shorttitle{Trannseries and
  summability}

\maketitle

\tableofcontents
\section{Introduction}

In the first part we give a rigorous and concise construction of a
space of transseries adequate for the study of a relatively large
class of ordinary differential and difference equations. The main
stages of the construction and the notations are similar to those in
\cite{EcalleNATO}, but the technical details rely on a generalized
contractivity principle which we prove for abstract multiseries
(Theorem~\ref{PFp}).  We then discuss how various classes of problems
are solved within transseries.  The last part of the paper is devoted
to a brief overview of generalized Borel summation, an isomorphism
used to associate actual functions to transseries, in the context of
ODEs, difference equations and the extension of some of these
techniques to PDEs.

Various spaces of transseries were introduced roughly at the same time
in logic (see \cite{VDD}), analysis \cite{Eca84}, \cite{EcalleNATO},
and the theory of surreal numbers \cite{Conway}; similar structures
were introduced independently by Berry and used as a powerful tool in
applied mathematics, see  \cite{Berry-hypB,BerryKeating}.

Informally, transseries are {\em asymptotic\footnote{An asymptotic
    expansion is one in which the terms are ordered decreasingly with
    respect to the order relation $\ll$, where $f_1(x)\ll f_2(x)$ if
    $f_1(x)=o(f_2(x))$}, finitely generated}\index{asymptotic
  expansion} combinations powers, logarithms and exponentials. It is a
remarkable fact that a wide class of functions can be asymptotically
described in terms of them. As Hardy noted \cite{Hardy} ``No function
has yet presented itself in analysis the laws of whose increase, in so
far as they can be stated at all, cannot be stated, so to say, in
logarithmic-exponential terms''. This universality of representation
can be seen as a byproduct of the fact that transseries are closed
under many common operations.

It is convenient to take $x^{-1}$, $x\to+\infty$, as the small
parameter in the transseries. When corresponding functions are
considered we change variables so that the limit is again $x\to
+\infty$.

A simple
example of a transseries generated by $x^{-1}$ and $e^{-x}$ is

$$\sum_{k,m=0}^{\infty}c_{km}e^{-kx}x^{-m}$$
where $c_{km}\in\CC$. This transseries is without logarithms or
log-free, and has {\em level }one since it contains no iterated
exponential; the simplest nontrivial transseries of level two is
$e^{e^x}$.  A more complicated example of a transseries
\index{transseries} of exponential level two, with level zero {\em
  generators} $x^{-1}$ and $x^{-\sqrt{2}}$, level one generator
$\exp(x)$, and level two generators $\exp(\sum_{k=0}^{\infty}c_k e^x
x^{-k})$ and $ \exp(-e^x)$ is

$$e^{\sum_{k=0}^{\infty}c_k e^x x^{-k}}+\sum_{k=0}^{\infty}
d_kx^{-k\sqrt{2}}+e^{-e^x}$$
Some examples of transseries-like expressions for $x\rightarrow
+\infty$, which are in fact {\em not} transseries, \index{transseries}
are $\sum_{k=0}^{\infty}x^k$ (it fails the asymptoticity condition)
and $\sum_{k=0}^{\!\infty}\!e^{-e^{nx}}$ (it does not have finitely
many generators, this property is described precisely in the sequel).

\subsection{Abstract multiseries}\index{transseries}
The underlying structure behind the condition of asymptoticity is that
of {\em well ordering}.  In order to formalize transseries
\index{transseries} and study their properties, it is useful to first
introduce and study more general abstract expansions, over a well
ordered set.

\subsubsection{Totally ordered sets; well ordered sets}
Let $A$ be an ordered set, with respect to $\le$. If $x\not\le y$ we write $x>y$ or $y<x$.
 $A$ is {\bf totally ordered} if any two elements are
comparable, i.e., if for any $x,y\in A$ we have $x\ge y$ or $y\ge
x$. If $A$ is not totally ordered, it is called
{\bf partially ordered}.

The set $A$ is {\bf well ordered} with respect to $>$ if every
  nonempty totally ordered subset ({\em chain}) of $A$ has a
  minimal element, i.e.

$$A'\subset A\implies \exists M\in A'\mbox{ such that } \forall
x\in A',\,\,M\le x$$

If any nonempty totally ordered subset of $A$ has a {\bf maximal} element, we say that
$A$ is well ordered with respect to $<$.
\subsubsection{Finite chain property} $A$ has the {\bf finite descending chain} property
if there is no {\em infinite strictly decreasing} sequence in $A$, in
other words if $f:\NN\mapsto A$ is decreasing, then $f$ is constant for
large $n$.

\begin{Proposition}\label{WOFC}
$A$ is well ordered with respect to ``$>$'' {\em iff} it has the finite
descending
chain property.
\end{Proposition}

\begin{proof}
  A strictly decreasing infinite sequence is obviously totally ordered
  and has no minimal element. For the converse, if there exists
  $A'\subset A$ such that $\forall\,x\in A'\,\exists\, y=:f(x)\in
  A',\,\,f(x)<x$ then for $x_0\in A'$, the sequence
  $\{f^{(n)}(x_0)\}_{n\in\NN}$ is an infinite descending chain in $A$.
\end{proof}

{\em Example: multi-indices.} $\NN$ is well ordered with respect to $>$,
and so is 

$$\NN^M-\bk_0:=\{\bk\in\ZZ^M:\bk\ge-\bk_0\}$$
with respect to the order relation $\bf m\ge n$ iff $m_i\ge
n_i\,\forall\,i\ge k$. Indeed an infinite descending sequence $\mathbf{n}_i$
would be infinitely descending on at least one component.

\begin{Proposition}\label{Va}
Let $\bk_0\in\ZZ^M$ be fixed. Any infinite set $A$ in $\NN^M-\bk_0$
contains a {\em strictly increasing} (infinite) sequence.
\end{Proposition}

\begin{proof}
  The set $A$ is unbounded, thus there must exist at least one component
  $i\le M$ so that the set $\{m_i:\mathbf{m}\in A\}$ is also unbounded;
  say $i=1$. We can then choose a sequence
  $S=\{\mathbf{m}_n\}_{n\in\NN}$ so that $(\mathbf{m}_n)_1$ is strictly
  increasing. If the set $\{(\mathbf{m}_n)_j: \mathbf{m}\in A;\,\,
  j>1\}$ is bounded, then there is a subsequence $S'$ of $S$ so that
  $(m_2,...,m_n)_{n'}$ is a constant vector. Then $S'$ is a strictly
  increasing sequence. Otherwise, one component, say
  $(\mathbf{m}_{n'})_2$ is unbounded, and we can choose a subsequence
  $S''$ so that $(m_1,m_2)_{n''}$ is increasing.  The argument continues
  in this fashion until in at most $M$ steps an increasing sequence is
  constructed.
\end{proof}
\begin{Corollary}
  \label{nononcomp} Any infinite set of multi-indices
  in $\NN^M$ contains at least two comparable elements.
\end{Corollary}

\begin{Corollary}
  \label{Lowerset} Let $A$ be a nonempty set of multi-indices in
  $\NN^M-{\bk_0}$.  There exists a unique and {\bf finite} minimizer set
  $\mathcal{M}_A$ such that none of its elements are comparable
and for any $a'\in A$ there is an $a\in \mathcal{M}_A$
  such that $a\le a'$.
\end{Corollary}
\begin{proof}
  Consider the set $C$ of all {\em maximal} totally ordered subsets of
  $A$ (every chain is contained in a maximal chain; also, in view of
  countability, Zorn's lemma is not needed). Let $\mathcal{M}_A$ be the set of the
  least elements of these chains, i.e. $\mathcal{M}_A=\{\min c:c\in C\}$.  Then
  $\mathcal{M}_A$ is finite. Indeed, otherwise, by Corollary~\ref{Lowerset} at
  least two elements in $\mathcal{M}_A$ such that $a'_1<a'_2$. But this
  contradicts the maximality of the chain whose least element was
  $a'_2$. It is clear that if $\mathcal{M}'_A$ is a minimizer then $\mathcal{M}'_A\supset
  \mathcal{M}_A$. Conversely if $m\in \mathcal{M}'_A\setminus \mathcal{M}_A$ then $m\nless a,\,\,\forall 
a\in \mathcal{M}_A$ contradicting the definition of $\mathcal{M}_A$.

\end{proof}
\subsubsection{Definition and properties of abstract multiseries}\label{SFst}
If $\mathcal{G}$ is a commutative group with an order relation, we
call it an {\bf abelian ordered group} if the order relation is
compatible with the group operation, i.e., $a\le A\mbox{ and }b\le
B\implies ab\le AB$ (e.g.  $\RR$ or $\ZZ^M$ with addition). Let
$\mathcal{G}$ be an abelian ordered group and let $\mu:\ZZ^M\mapsto
\mathcal{G}$ be a {\em decreasing group morphism}, i.e.,

(1) $\mu_{\mathbf{0}}=1$.

(2) $\mu_{\bk_1+\bk_2}=\mu_{\bk_1}\mu_{\bk_2}$.

(3) $\bk_1>\bk_2\implies \mu_{\bk_1}<\mu_{\bk_2}$.
Then $\mu(\ZZ)$ is the subgroup finitely generated by
$\mu_{\mathbf{e}_j};j=1,...,M$ where $\mathbf{e}_j$ is the
unit vector in the direction $j$ in $\ZZ^M$, and since
$\mathbf{e}_j>0$ it follows that

$$\mu_{\mathbf{e}_j}<1,\ j=1,...,M$$

In view of our final goal, a simple example to keep in mind is the
multiplicative group of {\em monomials} $\mathcal{G}_1$, generated by the functions $x^{-1/2}$, $x^{-1/3}$
and $e^{-x}$, for large positive $x$. The order relation on $\mathcal{G}_1$ is
$\mu_1<\mu_2$ if $|\mu_1(x)|<|\mu_2(x)|$ for all large $x$.  When $M=3$ 
we choose $\mu(\bk)=x^{-k_1/2-k_2/3}e^{-k_3 x}$.

\begin{Remark}\label{R0} The relation $\mu_{\bk_1}=\mu_{\bk_2}$
  induces an equivalence relation on $\ZZ^r$; we denote it by $\equiv$.
\end{Remark} 
For instance in $\mathcal{G}_1$, since $1/2$ and $1/3$ are not rationally
independent, there exist distinct $\bk',\bk$ so that
$\bk'\in\NN^3:\mu_{\bk'}=\mu_{\bk}$.

\begin{Remark}\label{finigen}
Clearly any choice of $\mu_i$ with $\mu_i<1$ for $i=1,...,M$ defines
an order preserving  morphism via $\mu(\bk)=\prod_{i=1}^M\mu_{i}^{k_i}$.
\end{Remark}

Ordered morphisms preserve well-ordering:
\begin{Proposition}\label{(4)}
Let $P\subset\ZZ^M$ be {\em well ordered} (an important
example for us is $P=\NN^M-\bk_0$)  and $\mu$ an order preserving
morphism. Then $\mu(P)$ is well ordered.
\end{Proposition}

\begin{proof}
  Assuming the contrary, let $J=\{\bk_n\}_{n\in\NN}$ be such that
  $\boldsymbol\mu_J=\{\boldsymbol\mu_{\bk_n}\}_{n\in\NN}$ is an infinite
  strictly {\em ascending} chain in $\mu(P)$. Then the index set $J$ is
  clearly infinite, and then, by Proposition~\ref{Va} it has a strictly
  increasing subsequence $J'$.  Then $\boldsymbol\mu_{J'}$ is a {\em
  descending} subsequence of $\boldsymbol\mu_J$, which is a
  contradiction.
\end{proof}
\begin{Corollary}\label{bilateral}
For any $\bk\in\ZZ^M$, the set $\{\bk'\in\NN^M-\bk_0:\mu_{\bk'}=\mu_\bk\}$
is finite. In particular, given $\bk''$, the set $\{\bk,\bk'\in\NN^M-\bk_0:\bk+\bk'=
\bk''$ is finite.
\end{Corollary}
\begin{proof}
  By Proposition~\ref{Va}, the contrary would imply the existence of
  an strictly increasing subsequence of $\bk'$, for which then
  $\{\mu_{\bk'}\}$ would be strictly decreasing, a contradiction. The
  last part follows if we take $\mu_{\bk}=\bk$.
\end{proof}

\begin{Proposition}\label{2.9}
  The space of formal series\footnote{i.e.  the space of real or complex
    {\em functions} on $\mu(\NN^M-\bk_0)$ with usual addition and
    convolution\index{convolution}  (\ref{innerprod1}).}
    $$\tilde{\mathcal{A}}(\mu_1,...,\mu_M)=
\tilde{\mathcal{A}}=\{\tilde{S}=\sum_{\bk\ge
    \bk_0}c_{\bk}\mu_{\bk}:
 \bk_0\in\ZZ^M,\,c_{\bk}\in\CC\}$$
is an algebra with respect to
componentwise multiplication by scalars, componentwise
addition, and the inner multiplication
\begin{equation}\label{innerprod}
  \tilde{S} \tilde{S}'=\sum_{\bk\ge \bk_0}\sum_{\bk'\ge
  \bk'_0}c_{\bk}c_{\bk'}\mu_{\bk+\bk'}=
\sum_{\bk''\ge \bk_0+\bk'_0}\mu_{\bk''}c_{\bk''}
\end{equation}
where
\begin{equation}\label{innerprod1}
c_{\bk''}=
\mathop{\sum_{\bk\ge \bk_0,\bk'\ge\bk_0}}
_{\bk+\bk'=\bk''}c_{\bk}c_{\bk'}
\end{equation}
The same is true for $\tilde{\mathcal{A}}(\mu_1,...,\mu_M)$ factored by
the equivalence relation
\begin{equation}
  \label{eq:defequi}
  \tilde{S}\equiv \tilde{S}'\iff \sum_{\mathbf{k}'\equiv\bk}(c_{\mathbf{k}'}-
c'_{\mathbf{k}'})=0\ \forall\,
\bk\ge \bk_0
\end{equation}
 replacing $\bk+\bk'=\bk''$ with $\bk+\bk'\equiv\bk''$ in
(\ref{innerprod}) (note that by Corollary~\ref{bilateral} the
equivalence classes have finitely many elements).
\end{Proposition}
\begin{proof}
 Straightforward.
\end{proof}
We define $\tilde{\mathcal{A}}_{\bk_0}$ the linear subspace 
of $\tilde{\mathcal{A}}$ for which $\bk_0$ is fixed.
\begin{Definition}
The sum
$$\tilde{S}_c=\sum_{\bk\ge \bk_0}c_{\bk}\mu_{\bk}$$
 is in {\em collected} form if, by definition, $c_\bk\ne 0\implies
\bk=\max\{\bk':\bk'\equiv\bk\}$, where the maximum is with respect to
the lexicographic order. (In other words the coefficients are collected
and assigned to the earliest $\mu$ in its equivalence class.)

It is then natural to represent the equivalence class (v.
(\ref{eq:defequi})) $\{\tilde{S}\}$ of $\tilde{S}$, in
$\tilde{\mathcal{A}}_{\bk_0}/\equiv$,  by $\tilde{S}_c$.
\end{Definition}
\begin{Corollary}\label{collected}
By Proposition \ref{2.9} every nonzero sum can be written in collected form.
\end{Corollary}
\subsection{Topology on multiseries}
A topology
 is introduced in the following way:
\begin{Definition}
  \label{Topology}
  {\rm The sequence $\tilde{S}^{(n)}$ in $\tilde{\mathcal{A}}_{\bk_0}$
    converges in the asymptotic topology \index{topology}if for any $\bk$,
    $c_{\bk}^{(n)}$ becomes  constant (with respect to $n$) eventually.} This
  induces a natural topology on $\tilde{\mathcal{A}}_{\bk_0}/\equiv$.
\end{Definition}
This topology \index{topology}is metrizable. Indeed, any bounded
$w:\NN^M\mapsto (0,\infty)$ such that $w(\bk)\rightarrow 0$
iff all $k_i\rightarrow \infty$ (e.g. $w(\bk)=\sum_{i\le
M}e^{-k_i}$) provides a translation-invariant distance
$$d(\tilde{S}^{(1)},\tilde{S}^{(2)})=\sup_{\bk\ge \bk_0}\theta(c_{\bk}^{(1)}-c_{\bk}^{(2)})w(\bk-\bk_0)$$ 
where $\theta(x)=0$ if $x= 0$
and is one otherwise.  A  Cauchy sequence in
$\tilde{\mathcal{A}}_{\bk_0}$ is clearly convergent, and in this sense
$\tilde{\mathcal{A}}_{\bk_0}$ is a complete topological algebra.
 From this point on, we assume $\mathcal{G}$ is
a {\bf totally ordered abelian group}. Let
$\tilde{S}\in\tilde{\mathcal{A}}$. 
\begin{Remark}
  A subgroup of $\mathcal{G}$ generated by $n$ elements $\mu_1<1,...,\mu_n<1$
is totally ordered and well ordered, and thus can be indexed by a set
of ordinals $\Omega$, in such a way that $\omega_1<\omega_2$
implies $\mu_{\omega_1}>\mu_{\omega_2}$. A sum
\begin{eqnarray}
  \label{eq:asympt-form}
  \tilde{S}=\sum_{\omega\in\Omega_{\tilde{S}}}c_{\omega}\mu_{\omega}
\end{eqnarray}
\z where we agree to omit from $\Omega_{\tilde{S}}$ all ordinals for which
$c_\omega=0$ is called the {\em asymptotic} form of $\tilde{S}$.
\end{Remark}
\begin{Definition} {\bf Dominant term.}
  Assume $\tilde{S}=\sum_{\bk\ge\bk_0} c_\bk
  \mu_\bk\in\mathcal{A}_{\bk_0}$ is presented in collected form.  The
  set $\mu_{\bk}:c_{\bk}\ne 0$ is then totally ordered and must have a
  maximal element $\mu_{\bk_1}$, by Proposition~\ref{(4)}.  We say
  that $c_{\bk_1}\mu_{\bk_1}=:\dom(\tilde{S})$ is the dominant term of
  $\tilde{S}$ and $\mu_{\bk_1}=:\mag(\tilde{S})$ is the {\em
    (dominant) magnitude} of $\tilde{S}$ (equivalently,
  $\mag(\tilde{S})=\mu_{\min\Omega_{\tilde{S}}}$). We allow for
  $\mag(\tilde{S})$ to be zero, {\bf iff} $\tilde{S}=0$.
\end{Definition}
 The following property is an immediate consequence of
Corollar\ref{collected}:
\begin{Remark}\label{R2.19}
For any nonzero $\tilde{S}$ we can write 
$$\tilde{S}=c_{\bk_1}\mag(\tilde{S})(1+\sum c'_{\bk'}\mu_{\bk'})
=\dom(\tilde{S})(1+\tilde{S}_1)$$
where all the terms in $\tilde{S}_1$ are less than one.
\end{Remark}
 \begin{Remark}\label{Rfinitude} The magnitude is continuous:
  if $\tilde{S}_{\bk}\in\tilde{\mathcal{A}}_{\bk_0}$  and
  $\tilde{S}_{\bk}\rightarrow\tilde{S}$ in the
  asymptotic topology, then
$\mag(\tilde{S}_{\bk})\rightarrow\mag(\tilde{S})$
(i.e.  $\exists\bk_1$ so that
  $\mag(\tilde{S})=\mag(\tilde{S}_{\bk}),\,\forall \bk\ge \bk_1$).
\end{Remark}
\begin{proof}
This follows immediately from the definition of the topology
and of $\mag(\cdot)$.
\end{proof}
The proposition below discusses the closure of
$\tilde{\mathcal{A}}_{\bk_0}$ under restricted {\em infinite sums}.
\begin{Proposition}
  \label{genclos}
 Let $\mathbf{j}_0,\mathbf{k}_0,\mathbf{l}_0\in\ZZ^M$ with
$\mathbf{k}_0+\mathbf{l}_0=\mathbf{j}_0$ and consider the sequence
in $\tilde{\mathcal{A}}_{\bk_0}$
$$\tilde{S}^{(\mathbf{m})}=\sum_{\bk\ge \bk_0}c_{\bk}^{(\mathbf{m})}\mu_{\bk}$$
and a fixed $T\in\tilde{\mathcal{A}}_{\mathbf{l}_0}$,
$$T=\sum_{\bk\ge \mathbf{l}_0}c'_{\bk}\mu_{\bk}$$
Then the sum (``blending'')
$$T(\tilde{S}):=\sum_{\bk\ge \mathbf{l}_0}c'_{\bk}\mu_{\bk}\cdot
\tilde{S}^{(\mathbf{k})}$$
obtained by replacing each $\mu_{\bk}$ in
$T$ by the product $\mu_{\bk}\cdot \tilde{S}^{(\mathbf{k})}$ is well
defined in $\tilde{\mathcal{A}}_{\bj_0}$ as the limit as
$\bk'\uparrow\infty$ of the convergent sequence of {\em truncates}
\begin{multline}
  T^{[\bk']}(\tilde{S})=\sum_{\mathbf{l}_0\le\bk\le\bk'}c'_{\bk}\mu_{\bk}\cdot
\tilde{S}^{(\mathbf{k})}=
\sum_{\mathbf{j}\ge\bj_0}\mu_{\mathbf{j}}\sum_{B^{(\bk')}_\mathbf{j}}c'_{\mathbf{m}_1}c_{\mathbf{m}_2}^{(\mathbf{m}_1)}=
\sum_{\mathbf{j}\ge\bj_0}d_{\mathbf{j}}^{(\bk')}\mu_{\mathbf{j}}
\end{multline}
where
$$B_{\mathbf{j}}^{(\bk')}=\{\mathbf{m}_1,\mathbf{m}_2:
\mathbf{m}_1+\mathbf{m}_2=\mathbf{j},
\mathbf{l}_0\le\mathbf{m}_1\le\bk', \ \mathbf{m}_2\ge\bk_0\}
$$
\end{Proposition}
\begin{proof}\label{p} 
Given $\mathbf{j}$, the coefficient $d_{\mathbf{j}}^{(\bk')}$ is
constant for large  $\bk'$. Indeed, in the expression of
$B_{\mathbf{j}}^{(\bk')}$ we have $\mathbf{l}_0\le\mathbf{m}_2=
\mathbf{j}-\mathbf{m}_1\le \mathbf{j}-\mathbf{k}_0$ and 
similarly $\bk_0\le\mathbf{m}_2\le \mathbf{j}-\mathbf{l}_0$ and
therefore there is a bound independent of $\bk'$ on the number of
elements in the set $B^{(\bk')}_{\mathbf{j}}$. On the other hand, we
obviously have $B^{(\bk')}_{\mathbf{j}}\subset
B^{(\bk'')}_{\mathbf{j}}$ if $\bk''>\bk'$. Thus the set
$B^{(\bk')}_{\mathbf{j}}$ is constant if $\bk'$ is large enough, and
thus $d_{\mathbf{j}}^{(\bk')}$ is constant for all large $\bk'$, which
means the sums $T^{[\bk']}(\tilde{S})$ are convergent in the
asymptotic topology.
\end{proof}
{\bf Note}. The condition that $\mag(\tilde{S})^{(m)}$
 decreases strictly in $m$ does {\bf not} suffice
for $\sum_{m\ge 0}c_m\tilde{S}^{(m)}$ to be well defined. 
Indeed, the terms $\tilde{S}^{(m)}=x^{-m}+e^{-x}\in \mathcal{G}_1$ 
(cf. \S~\ref{SFst}) have strictly decreasing
magnitudes and yet the formal expression $\sum_{m\ge 0}(x^{-m}+e^{-x})$
is meaningless.
\subsection{Contractive operators}\label{SubContr}
\begin{Definition}\label{DFp}
Let $J$ be a linear operator from $\tilde{\mathcal{A}}_{\bk_0}$ or from
one of its subspaces, to $\tilde{\mathcal{A}}_{\bk_0}$,
\begin{equation}
  \label{EJlinear}
  J\tilde{S} =J\sum_{\bk\ge\bk_0}c_{\bk}\mu_{\bk}=\sum_{\bk\ge\bk_0}c_{\bk}J\mu_{\bk}
\end{equation}
 Then $J$
is called asymptotically contractive\index{asymptotically contractive} on $\tilde{\mathcal{A}}_{\bk_0}$ if
\begin{equation}
  \label{Econtractformula}
 J\mu_{\mathbf{j}}
=
\sum_{\mathbf{p}>0}c_{\mathbf{j;p}}\mu_{\bj+\mathbf{p}}
\end{equation}
thus
\begin{equation}
  \label{Econtractformula1}
 J\tilde{S}=\sum_{\bk\ge\bk_0}\mu_{\bk}
\mathop{\sum_{\bj+\mathbf{p}\equiv\bk}}_{\mathbf{p}>0;\bj\ge\bk_0}c_{\bj;\mathbf{p}}c_{\bj}
\end{equation}
\end{Definition}
\z We note that by (\ref{Econtractformula}) and
Proposition~\ref{genclos}, $J$ is well defined.
\begin{Definition}
  \label{nonlinear.contrac}
  The linear or nonlinear operator $J$ is (asymptotically)
contractive\index{asymptotically contractive}
  in the set $A\subset\mathcal{A}_{\bk_0}$ if $J:A\mapsto A$ and the
  following condition holds. Let $f_1$ and $f_2$ in $A$ be arbitrary and let
\begin{eqnarray}
    \label{Econtr.notat}
 f_1-f_2=\sum_{\bk\ge \bk_0}c_{\bk}\mu_{\bk}
  \end{eqnarray}
 Then
  \begin{eqnarray}
    \label{Econtr.nonl}
    J(f_1)-J(f_2)=\sum_{\bk\ge \bk_0}c_{\bk}\mu_{\bk}\tilde{S}_{\bk}
  \end{eqnarray}
 where $\mag(\tilde{S}_{\bk})=\mu_{\mathbf{p_\bk}}$ for some $\mathbf{p_\bk}>0$.\end{Definition}
\begin{Remark}\label{Rlinearity}
 The sum of asymptotically contractive\index{asymptotically contractive} operators is contractive; the
  composition
of contractive operators, whenever defined, is contractive.
\end{Remark}
\begin{Theorem}\label{PFp}
(i) If $J$ is linear and contractive\index{asymptotically contractive} on $\tilde{\mathcal{A}}_{\bk_0}$
  then for any $\tilde{S}_0\in\tilde{\mathcal{A}}_{\bk_0}$ the fixed
  point equation $\tilde{S}=J\tilde{S}+\tilde{S}_0$ has a unique
  solution $\tilde{S}\in\tilde{\mathcal{A}}_{\bk_0}$. 
  
  (ii) In general, if $A\in\mathcal{A}_{\bk_0}$ is closed and
  $J:A\mapsto A$ is a (linear or nonlinear) contractive\index{asymptotically contractive} operator on $A$,
  then $\tilde{S}=J(\tilde{S})$ has a unique solution is $A$.
\end{Theorem}
\begin{proof} (i) 
  Uniqueness: if we had two solutions $\tilde{S}_1$ and $\tilde{S}_2$ we
  would get
  $\mag(\tilde{S}_1-\tilde{S}_2)=\mag(J(\tilde{S}_1-\tilde{S}_2))<\mag(\tilde{S}_1-\tilde{S}_2)$
  by (\ref{Econtractformula}). We show that $J^n\tilde{S}_0\rightarrow
  0$ in the asymptotic topology implying that $\sum_{n=0}^{\infty}J^n
  \tilde{S}_0$ is convergent. We have
  \begin{align}\label{EJn(delta)}
    J^n\tilde{S}=\sum_{\bk\ge\bk_0}\mu_{\bk}
 \mathop{\sum_{\mathbf{p_1}>0,...,\mathbf{p}_n>0;\bj\ge\bk_0}}_{\bj+\mathbf{p}_1+\cdots+\mathbf{p}_n\equiv\bk}const_{\mathbf{p_1};...,\mathbf{p_n};\mathbf{j}}
  \end{align}
  and since $\left|\mathbf{p}_1+\cdots+\mathbf{p}_n\right|>n$ and since the set
 $\{\bk':\bk'\equiv\bk\}$ is finite by Corollary~\ref{bilateral}, for
 each $\bk$ the condition $\bj+\mathbf{p}_1+\cdots+\mathbf{p}_n\equiv\bk
 $ becomes impossible if $n$ exceeds some $n_0$, and then the
 coefficient of $\mu_\bk$ in $J^n\tilde{S}$ is zero for $n>n_0$.
 
(ii) Uniqueness follows in the same way
 as in the linear case. For existence, prove the convergence of the
 recurrence $f_{n+1}=J(f_n)$.
 
 With
 $f_{n}-f_{n-1}=\delta_{n-1}=\sum_{\bk\ge\bk_0}c_\bk^{(n-1)}\mu_\bk$
 we have, for some coefficients $C_{\mathbf{m};\mathbf{m}_1}^{(n-1)}$

\begin{multline}
  \label{En->n+1}
\delta_n=J(f_{n-1}+\delta_{n-1})-J(f_{n-1})=
\sum_{\bk\ge\bk_0}c_\bk^{(n-1)}\mu_\bk\tilde{S}^{(n-1)}_\bk\\=
\sum_{\bk\ge\bk_0}\mu_\bk\mathop{\sum_{\mathbf{m}+\mathbf{p}=\bk}}_{\mathbf{m}\ge\bk_0;\mathbf{p}>0}c_{\mathbf{m}}^{(n-1)}
C_{\mathbf{m};\mathbf{p}}^{(n-1)}:=\sum_{\bk\ge\bk_0}c_\bk^{(n)}\mu_\bk
\end{multline}

\z and therefore $\delta_n$ has an expression similar
to (\ref{EJn(delta)}),

$$\delta_n=\sum_{\bk\ge\bk_0}\mu_{\bk}
 \mathop{\sum_{\mathbf{p_1}>0,...,\mathbf{p}_n>0;\bj\ge\bk_0}}_{\bj+\mathbf{p}_1+\cdots+\mathbf{p}_n\equiv\bk}const_{\mathbf{p_1};...,\mathbf{p_n};\mathbf{j}}$$

\z Consequently $\delta_n\rightarrow 0$ and, as before,
it follows that $\sum_n\delta_n$ converges. 
\end{proof}
\begin{Corollary}\label{compo_n}
Let $\tilde{S}\in\mathcal{A}_0$ be arbitrary and
$\tilde{S}_n=\sum_{\bk>0}c_{\bk;n}\mu_{\bk}\in\mathcal{A}_0$ for
$n\in\NN$. Then the operator defined by

$$J(y)= \tilde{S}+\sum_{n\ge 2}\tilde{S}_{n-2}y^{n}$$
  \z is contractive\index{asymptotically contractive} in the set $\{y:\mag(y)<1\}$.
\end{Corollary}

\begin{proof}
  We have 

$$J(y+\delta)-J(y)=\delta\sum_{n\ge 2}\tilde{S}_{n-2}\left(
\sum_{j=1}^{n-1}y^j\delta^{n-j}\right)$$

\end{proof}

\subsubsection{The field of finitely generated formal series}

Let $\mathcal{G}$ be a totally ordered abelian group. We now define the algebra:

\begin{equation}
  \label{eq:defSapprox}
   \stackrel{\approx}{S}=\mathop
{\mathop{\bigcup_{\bk_0\in\ZZ^M}}_{\mu_1<1,...,\mu_M<1}}_{M\in\NN}
\tilde{\mathcal{A}}_{\bk_0}(\mu_1,...,\mu_M)
\end{equation}

\z modulo the obvious inclusions, and with the induced topology
(convergence in $\stackrel{\approx}{S}$ means convergence in one of the
$\tilde{\mathcal{A}}_{\bk_0}(\mu_1,...,\mu_M)$).

{\em Product form}. 
\begin{Proposition}\label{Pstd}
  Any $\tilde{S}\in \stackrel{\approx}{S}$ can be written in the form
  \begin{equation}
    \label{Estd}
    c\,\mag(\tilde{S})\left(1+\sum_{\bk> 0}c_{\bk}\mu_{\bk}\right)
  \end{equation}
  
  \z i.e., 
\begin{equation}
    \label{Estd2}
    c\,\mag(\tilde{S})(1+\tilde{S}_1)
  \end{equation}
  
\z where $\tilde{S}_1\in\tilde{\mathcal{A}}_{\bk_0}$ for some
$\bk_0> 0$ (cf. also Remark \ref{R2.19}) and $\mag(\tilde{S}_1)<1$. 
\end{Proposition}
\begin{proof}
  We have, by Remark \ref{R2.19},
 \begin{align}\label{caseI}
  \tilde{S}=
c\,\mag(\tilde{S})\left(1+\sum_{\bk\ge \bk_0}c'_{\bk}\mu_1^{k_1}\cdots
\mu_M^{k_M}\mag(\tilde{S})^{-1}\right)
\end{align}

\z where all the elements in the last sum are less than one.

  Let $A$ be the set of multi-indices in the sum in (\ref{caseI}) for which {\em some} $k_i <0$ and let $A'$   be
its minimizer in the sense of Corollary~\ref{Lowerset}, a finite set. We now
consider the extended set of generators

$$\{\overline{\mu}_i:i\le M'\}:=\{\nu\mag(\tilde{S})^{-1}:\nu=\mu_i\mbox{ with }i\le M\mbox{ or
  }\nu=\mu_\bk\mbox{ with }\bk\in A'\}$$

We clearly have $\overline{\mu}_i<1$. By the definition of $A'$, for any
term in the sum in (\ref{caseI}) either $\bk>0$ or
else $\bk=\bk'+\bk''$ with $\bk'\in A'$ and $\bk''\ge 0$. In both cases we
have  $c_{\bk}\mu_{\bk}=c_{\bk'}\overline{\mu}_{\bk'}$ with
$\bk'\in\ZZ^{M'}$ and $\bk'>0$. Thus $\tilde{S}$ can be rewritten in the
form

$$c\,\mag(\tilde{S})\left(1+\sum_{\bk> 0}c_{\bk}\mu_{\bk}\right)\ \
(\bk\in\NN^{M'})$$

\z where the assumptions of the Proposition are satisfied.

\end{proof}
\begin{Proposition}
  \label{Field}
 $\stackrel{\approx}{S}$ is a field.
\end{Proposition}

\begin{proof}
  The only property that needs verification is the existence of a
  reciprocal for any nonzero $\tilde{S}$. Using Proposition~\ref{Pstd}
  we only need to consider the case when
$$\tilde{S}=1+\sum_{\bk>0}c_{\bk}\mu^\bk$$

\z Since multiplication by $t$ is manifestly contractive\index{asymptotically contractive} (see
\S~\ref{SubContr}), $\tilde{S}^{-1}$ is the solution (unique by
Proposition~\ref{PFp}) of

$$\tilde{S}^{-1}=1-t\tilde{S}^{-1}$$

\end{proof}

{\em Closure under infinite sums.} 
\begin{Corollary}
  \label{infisum} (i) Let $\bk_0>0$ and
$\tilde{S}\in\tilde{\mathcal{A}}_{\bk_0}$ and $\{c_n\}_{n\in\NN}\in\CC$
be any sequence. Then

$$\sum_{n=0}^{\infty} c_n \tilde{S}^n\in\tilde{\mathcal{A}}_{\bk_0}$$

(ii) More generally,
if $\tilde{S}_{01},...,\tilde{S}_{0M}$ are of the form $\tilde{S}_0$ and
$\{c_{\bk}\}_{\bk\in\NN^M}$ is a multi-sequence of constants, then
$\sum_{\bk\ge 0}c_\bk \tilde{S}_0^{\bk}=\sum_{\bk\ge 0}c_\bk
\tilde{S}_{01}^{k_1}\cdots \tilde{S}_{0M}^{k_M}$ is well defined.
\end{Corollary}

\begin{proof}
For some $C\in\CC$ we have $\tilde{S}=C\mag(\tilde{S})(1+t)=C(1+t)\mu$.
Since $\mu<1$ we have $\sum_{n=0}^{\infty} c_n
C^n\mu^n\in\tilde{\mathcal{A}}_0(\mu) $  and by Proposition~\ref{genclos}
$$\sum_{n=0}^{\infty} c_n \tilde{S}^n=
\sum_{n=0}^{\infty} c_n C^n\mu_{n\bk_1}(1+t)^n\in\stackrel{\approx}{S}$$
\end{proof}

{\em Formal series with real coefficients. Order relation}. Let
$\stackrel{\approx}{S}_\RR$ the subfield consisting in $\tilde{S}\in
\stackrel{\approx}{S}_\RR$ which have real
coefficients.  We say that $\stackrel{\approx}{S}_\RR\ni \tilde{S}>0$ if
$\dom(\tilde{S})/\mag(\tilde{S})>0$. Then every nonzero $\tilde{S}\in \stackrel{\approx}{S}_\RR$ is
either positive or else $-\tilde{S}$ is positive. This induces a {\em total}
order relation on $\stackrel{\approx}{S}_\RR$, by writing $\tilde{S}_1>\tilde{S}_2$ if
$\tilde{S}_1-\tilde{S}_2>0$.

\vfill\eject
\subsection{Inductive construction of logarithm-free transseries \index{transseries}}\label{Slogfree}
\subsubsection{Transseries}

Transseries are constructed as a special instance of abstract series
in which the ablelian ordered group is constructed inductively.

In constructing spaces of transseries, one aims at constructing
differential fields containing $x^{-1}$, closed under all operations
of importance for a certain class of problems operations. Smaller
closed spaces can be endowed with better overall properties.

The construction presented below differs in a number of technical
respects from the one of \'Ecalle, and the transseries space
constructed here is smaller than his. Still some of the construction
steps and the structure of the final object are similar enough to
\'Ecalle's, to justify using his terminology and notations.

\subsubsection{\'Ecalle's notation} 

\begin{itemize}
\item $\smallmo$ ---small transmonomial.
\item  $\largemo$ ---large transmonomial.
\item $\anymo$ ---any transmonomial, large or small.
\item
$\smalltr$ ---small transseries. \index{transseries}
\item
$\largetr$ ---large transseries. \index{transseries}
\item
$\anytr$
---any transseries, \index{transseries} small or large.

\end{itemize}

\subsubsection{Level 0: power series}\label{Level0} Let $x$  be large and positive
and let $\mathcal{G}$ be the totally ordered multiplicative group
$(x^{\sigma},\,\cdot, \,\ll), \sigma\in\RR$, with $x^{\sigma_1}\ll
x^{\sigma_2}$ if $x^{\sigma_1} =o(x^{\sigma_2})$ as
$x\rightarrow\infty$, i.e., if $\sigma_1<\sigma_2$. The space of level
zero log-free transseries \index{transseries} is by definition
$\tilde{\mathcal{T}}^{[0]}=\stackrel{\approx}{S}(\mathcal{G})$.  By
Proposition~\ref{Field}, $\tilde{\mathcal{T}}^{[0]}$ is a field.

If  $\tilde{T}\in\tilde{\mathcal{T}}^{[0]}$, then $\tilde{T}=\anymo$
iff $\tilde{T}=x^{\sigma}$ for some $\sigma\ne 0$, $\tilde{T}=\largemo$
if $\sigma>0$ and $\tilde{T}=\smallmo$ if $\sigma<0$.

The general element of $\tilde{\mathcal{T}}^{[0]}$ is a level zero
transseries \index{transseries}, $\anytr^{[0]}$ or $\anytr$ in short.  We have

\begin{equation}
  \label{defanytr}
  \anytr=\sum_{\bk\ge \bk_0}c_{\bk}\smallmo^{\bk}
\end{equation}

There are  two order
relations: $<$ and $\ll$ on $\tilde{T}\in\tilde{\mathcal{T}}^{[0]}$. 
We have $\anytr_1\ll \anytr_2$ iff $\mag(\anytr_1)\ll \mag(\anytr_2)$
(the sign of the leading coefficient is immaterial)
and $\anytr> 0$ if ($\anytr\ne 0$ and) the real number 
$\dom(\anytr)/\mag(\anytr)$ is positive.

\begin{Definition}
  \label{smalllarge}
{\rm A transseries \index{transseries} is {\bf small}, i.e. $\anytr=\smalltr$ iff in
(\ref{defanytr}) we have $c_{\bk}=0$ whenever $\smallmo^{\bk}\not\ll 1$.
Correspondingly, transseries \index{transseries} is {\bf large}, i.e. $\anytr=\largetr$ iff
in (\ref{defanytr}) we have $c_{\bk}=0$ whenever $\smallmo^{\bk}\not\gg
1$. We note that $\anytr=\smalltr$ iff $\mag(\anytr)\ll 1$ (there is an
asymmetry: the condition $\mag(\anytr)\gg 1$ does {\em not} imply
$\anytr=\largetr$, since it does not prevent the presence of small
terms in $\anytr$). } Any transseries \index{transseries} can then be written uniquely as

\begin{multline}
  \label{Edecsmalllarge}
  \anytr=\sum_{\bk\ge\bk_0}c_{\bk}\smallmo^{\bk}
=\sum_{\bk\ge\bk_0;\ \smallmo^{\bk}>1}c_{\bk}\smallmo^{\bk}
+const+\sum_{\bk\ge\bk_0;\ \smallmo^{\bk}<1}c_{\bk}\smallmo^{\bk}\\
=\largetr+const+\smalltr:=L(\anytr)+C(\anytr)+s(\anytr)
\end{multline}
\end{Definition}

\subsubsection{Level 1: Exponential power series}\label{Level1} The set $\mathcal{G}^{[1]}$
of {transmonomials of exponentiality one} consists by definition in the
formal expressions

$$\anymo^{[1]}=\anymo^{[0]}\exp({\largetr^{[0]}}),\ \
\anymo^{[0]},\largetr^{[0]}\in \tilde{\mathcal{T}}^{[0]}$$

\z where we allow for $\largetr^{[0]}=0$ and set $\exp(0)=1$. With
respect to the operation

$$\anymo_1^{[0]}\exp(\largetr_1^{[0]})
\anymo_2^{[0]}\exp(\largetr_2^{[0]})=
(\anymo_1^{[0]}\anymo_2^{[0]})\exp(\largetr_1^{[0]}+\largetr_2^{[0]})$$
\z we see that $\mathcal{G}^{[1]}$ is a commutative group. 

The order relations are introduced in the following way.

\begin{multline}\label{Eorderrel}
  \anymo_1\exp(\largetr_1)\gg
  \anymo_2\exp(\largetr_2)\\ \mbox{ iff }
  \left(\largetr_1>\largetr_2\right)\mbox{ or }
  \left(\largetr_1= \largetr_2\mbox{ and
      }\anymo_1\gg\anymo_2\right)
\end{multline}

\z In particular, if $\largetr_1^{[0]}$ is positive, then
$\anymo_1^{[0]}\exp(\largetr_1^{[0]})\gg 1$.

The second order relation, $>$, is defined by 

$$\anymo_1^{[0]}\exp(\largetr_1^{[0]})>0\iff\anymo_1^{[0]}>0$$

It is straightforward to check that $(\mathcal{G}^{[1]},\cdot,\gg)$
is an abelian ordered group. The abelian ordered
group of zero level monomials,  $(\mathcal{G}^{[0]},\cdot,\gg)$,
is naturally identified with the set of transmonomials for which
$\largetr^{[0]}=0$.

The space $\tilde{\mathcal{T}}^{[1]}$ of level one transseries
\index{transseries} is by definition
$\stackrel{\approx}{S}(\mathcal{G}^{[1]})$.  By Proposition~\ref{Field},
$\tilde{\mathcal{T}}^{[1]}$ is a field.  By construction, the space
$\tilde{\mathcal{T}}^{[0]}$ is embedded in
$\tilde{\mathcal{T}}^{[1]}$. Formula (\ref{defanytr}) is the general
expression of a level one transseries \index{transseries}, where now
$\smallmo$ is a transmonomial of level one.  The two order relations on
transseries \index{transseries} are the ones induced by transmonomials,
namely

\begin{equation}
  \label{Eord2}
  \anytr\gg 1\iff \mag(\anytr)\gg 1  \ \mbox{ and }\ \anytr> 0 \iff \dom(\anytr)/\mag(\anytr)>0
\end{equation}

\subsubsection{Induction step: level $n$ transseries \index{transseries}} Assuming the
transseries \index{transseries} of level $\le n-1$ are constructed, transseries \index{transseries} of level $n$
together with the order relation, are constructed exactly as in
\S~\ref{Level1}, replacing $[0]$ by $[n-1]$ and $[1]$ by $[n]$. 
The group $\mathcal{G}^{[1]}$ of transmonomials of order at most $n$
consists in expressions of the form

\begin{equation}
  \label{Ecannon}
  \anymo^{[n]}=x^{\sigma}\exp(\largetr^{[n-1]})
\end{equation}

\z where $\largetr^{[n-1]}$ is either zero or a large transseries of
level $n-1$ with the multiplication:

\begin{equation}
  \label{multcannon}
  x^{\sigma_1}\exp(\largetr_1^{[n-1]})x^{\sigma_2}\exp(\largetr_2^{[n-1]})
=x^{\sigma_1+\sigma_2}\exp(\largetr_1^{[n-1]}+\largetr_2^{[n-1]})
\end{equation}

\z The order relation is given by

\begin{align}
  \label{ineqcannon}
  x^{\sigma_1}\exp(\largetr_1^{[n-1]})\gg
x^{\sigma_2}\exp(\largetr_2^{[n-1]})
\iff \\
\Big(\largetr_1^{[n-1]}>\largetr_2^{[n-1]}\Big) \ \mbox{or} \ 
\Big(\largetr_1^{[n-1]}=\largetr_2^{[n-1]} \ \mbox{and}\ \sigma_1>\sigma_2\Big) 
\end{align}

\begin{equation}
  \label{defanytrn}
  \anytr^{[n]}=\sum_{\bk\ge \bk_0}c_{\bk}(\smallmo^{[n]})^{\bk}
\end{equation}

 As
in \S~\ref{Level1}, $\tilde{\mathcal{T}}^{[n-1]}$ is naturally embedded in
$\tilde{\mathcal{T}}^{[n]}$.

\subsubsection{General log-free transseries \index{transseries}, 
  $\stackrel{\approx}{\mathcal{T}}$} This is the space of arbitrary
level transseries, the inductive limit of the finite level spaces of
transseries \index{transseries}:

$$\stackrel{\approx}{\mathcal{T}}=\bigcup_{n=0}^{\infty}
\tilde{\mathcal{T}}^{[n]}$$

\z Clearly $\stackrel{\approx}{\mathcal{T}}$ is a field.  The order
relation is the one inherited from $\tilde{\mathcal{T}}^{[n]}$. The
topology \index{topology}is also that of an inductive limit, namely a sequence
converges iff it converges in $\tilde{\mathcal{T}}^{[n]}$ for {\em
  some} $n$.

\subsubsection{Further properties of transseries}

\z {\em Definition}. The level $l(\anytr)$ of $\anytr$ is $n$ if
$\anytr\in\tilde{\mathcal{T}}^{[n]}$ and
$\anytr\not\in\tilde{\mathcal{T}}^{[n-1]}$.

\begin{Proposition}
  \label{Rform0}
  
  If $n=l(\largetr_1)>l(\largetr_2)$  then
  $\largetr_1\gg\largetr_2$.

\end{Proposition}
\begin{proof}
 We may clearly take $n\ge 1$.  Since (by definition)
$\largetr\gg 1$ we must have, in
 particular,
$\dom(\largetr)=cx^{\sigma}\exp(\largetr')$ with 
$\largetr'\ge 0$. 
By induction, and the assumption $l(\largetr_1)=n$ we must have 
$\largetr'_1>0$ and 
$l(\largetr_1')=n-1$. The proposition follows since,
by again by the induction step, $\largetr_1'\gg\largetr_2'$.
\end{proof}
\begin{Remark}
  \label{Rform}

If $\anymo$ is of level no less than $1$, then either $\anymo$ is large, and
then 
$\anymo\gg x^\alpha,\,\forall\alpha\in\RR$ or else  $\anymo$ is small, and then
$\anymo\ll x^{-\alpha},\,\forall\alpha\in\RR$.

\end{Remark}

\begin{Remark}
  \label{Rfinige1}
  
  We can define  {\em generating monomials} of $0\ne \anytr\in\tilde{\mathcal{T}}^{[n]}$  a minimal
  subgroup $\mathcal{G}=\mathcal{G}(\anytr)$ of ${\mathcal{G}}^{[n]}$
  with the following properties: 
\begin{itemize} \item $\anytr\in
  \stackrel{\approx}{S}(\mathcal{G})$; \item $x^\sigma_1
  \exp(\largetr_1)\in \mathcal{G}$ implies $x^\sigma_1 \in \mathcal{G}$
  and, if $\largetr_1\ne 0$, then $\mathcal{G}\supset \mathcal{G}(\largetr_1)$.
\end{itemize}

\z By induction we see that $\mathcal{G}(\anytr)$ is finitely generated
for any $\anytr\in \mathcal{T}^{[n]}$.

\end{Remark}

\subsubsection{Closure of $\stackrel{\approx}{\mathcal{T}}$ 
  under composition and differentiation} 
\begin{Proposition}
  \label{Pdifferentiation}
$\stackrel{\approx}{\mathcal{T}}$ and $\mathcal{T}^{[n]};\ n\in\NN$ are differential fields.
\end{Proposition}

\begin{proof}
Differentiation
${\mathcal{D}}=\frac{\mathrm{d}}{\mathrm{d}x}$
is introduced
inductively on $\stackrel{\approx}{\mathcal{T}}$, as term by term
differentiation, in the following way. Differentiation in
$\tilde{\mathcal{T}}^{[0]}$ is defined as:

\begin{equation}
  \label{defdif0}
  {\mathcal{D}}\anytr=\sum_{\bk\ge \bk_0}c_{\bk}{\mathcal{D}}\smallmo^{\bk}
\end{equation}

\z where, as mentioned in \S\ref{Level0} we have
$\smallmo=x^{-\sigma}$ for some $\sigma\in\RR^+$ and, in the natural
way, we set ${\mathcal{D}} x^{-\sigma}=-\sigma x^{-\sigma-1}$. This makes
${\mathcal{D}}\anytr\in \tilde{\mathcal{T}}^{[0]}$, and the generating
transmonomials of 
${\mathcal{D}} \anytr$ are those of $\anytr$ together with $x^{-1}$.

We assume by induction that differentiation
${\mathcal{D}}:\tilde{\mathcal{T}}^{[n-1]}\mapsto
\tilde{\mathcal{T}}^{[n-1]}$ has been defined for all transseries \index{transseries} of
level at most $n-1$. (In particular,
${\mathcal{D}}\anytr$ is finitely generated.) We define 

\begin{multline}
  \label{def:deriv:n}
  {\mathcal{D}}
  \left(\anymo^{[n]}\right)={\mathcal{D}}\left(x^{\sigma}\exp(\largetr^{[n-1]})\right)=
\sigma x^{\sigma-1}\exp(\largetr^{[n-1]})
\\+x^{\sigma}{\mathcal{D}}\largetr^{[n-1]}\exp(\largetr^{[n-1]})
\end{multline}

\z A level $n$ transseries \index{transseries} is 

\begin{equation}
  \label{defanytrn2}
  \anytr=\sum_{\bk\ge \bk_0}c_{\bk}\smallmo^{\bk}
=\sum_{\bk\ge \bk_0}c_{\bk}\prod_{j=1}^{M}\smallmo_j^{k_j}
\end{equation}

\z and we write in a natural way
\begin{multline}
  \label{defanytrnp}
  {\mathcal{D}}\anytr^{[n]}
=\sum_{\bk\ge \bk_0}c_{\bk} \sum_{m=1}^M k_m\smallmo_m^{k_m-1}
{\mathcal{D}} \smallmo_m
\prod_{m\ne j=1}^{M}\smallmo_j^{k_j}\\=
\sum_{m=1}^M\smallmo_m^{-1}\sum_{\bk\ge \bk_0}c_{\bk}  k_m\smallmo_m^{k_m-1}
{\mathcal{D}} \smallmo_m
\prod_{m\ne j=1}^{M}\smallmo_j^{k_j}
\end{multline}
\z and the result follows from the induction hypothesis, since

\begin{equation}
  \label{oneterm}
  \sum_{\bk\ge \bk_0}c_{\bk}  k_m \smallmo^{\bk}
\frac{{\mathcal{D}} \smallmo_m}{ \smallmo_m}=\frac{{\mathcal{D}} \smallmo_m}{\smallmo_m}\sum_{\bk\ge \bk_0}c_{\bk}  k_m \smallmo^{\bk}
\in\tilde{\mathcal{T}}^{[n]}
\end{equation}

\end{proof}
\begin{Corollary}
  \label{Cfinigen} If $\mathcal{G}_{\anytr}$ is the group (finitely) generated by
all generators in any of the levels of $\anytr$, then ${\mathcal{D}}\anytr$
is generated by the transmonomials of $\mathcal{G}_{\anytr}$ together possibly
with $x^{-1}$. If $\anytr\ne \mathrm{Const}.$ then $l(\anytr)=l(\anytr')$.
\end{Corollary}
\begin{proof}
  Immediate induction; cf. also the beginning of the proof of
  Proposition
~\ref{Pdifferentiation}.
\end{proof}
The properties of differentiation are the usual ones:
\begin{Proposition}
  \label{Pdiffprop}
${\mathcal{D}}(fg)=g{\mathcal{D}} f+f{\mathcal{D}} g$, ${\mathcal{D}} const=0$
and ${\mathcal{D}}(f\circ g)=({\mathcal{D}} f)\circ g{\mathcal{D}} g$
(for composition, see \S~\ref{Scompo}).
\end{Proposition}
\begin{proof}
  The proof is straightforward induction.
\end{proof}

In the space of transseries, \index{transseries} differentiation is also
compatible with the order relation, a property which is not true in
general, in function spaces. 

\begin{Proposition}
  \label{Porrder} For any $\largetr_i, i=1,2$ and $\smalltr_i, i=1,2$ 
we have
  $$\largetr_1\gg\largetr_2 \Leftrightarrow \largetr_1'\gg\largetr_2' $$
$$\smalltr_1\gg\smalltr_2 \Leftrightarrow \smalltr_1'\gg\smalltr_2' $$
$$\largetr_1'\gg\smalltr_1'$$ 

\end{Proposition}
\begin{proof}
  The proof is by induction. It is true for power series, which are the
  level zero transseries \index{transseries}. Assume the property holds for transseries \index{transseries} of
  level $\le n-1$, and first prove the result for {\em transmonomials}
  of order $n$, i.e. for the case when for $i=1,2$
  $\anytr_{i}=\anymo_{i}^{[n]}=x^{\sigma_{i}}\exp(\largetr_{i})$, where at least one of $\largetr_{i}$ has  level $n-1$.  

 We have to evaluate

$$\exp(\largetr_1-\largetr_2)x^{\sigma_1-\sigma_2}
\frac{\largetr_1'+\sigma_1x^{-1}}{\largetr_2'+\sigma_2x^{-1}}$$
and it is plain that we can assume without loss of generality
that $\largetr_1>0$.

(1) If $l(\largetr_1-\largetr_2)=n-1$ then $M\gg 1$ by
Proposition~\ref{Rform0}. The remaining case is that for $i=1,2$ we have
$l(\largetr_{i})=n-1$ but $\largetr_{i}$ are equal  or else $l(\largetr_1-\largetr_2)\le n-2$
(which obviously requires $n\ge 2$). Let
$\Delta=\largetr_1-\largetr_2$. We have
$l(\Delta)<l(\largetr_1)$ and also $\largetr_1\gg x^{\alpha}$
for some $\alpha>0$, thus by Proposition~\ref{Rform0} and the induction
hypothesis we have
$\largetr_1' \gg \Delta'+\sigma x^{-1}$
$$\largetr_1'+\Delta'+\frac{\sigma_2}{x}=\largetr_1'(1+\smalltr)$$
\z  thus

$$\exp(\largetr_1-\largetr_2)x^{\sigma_1-\sigma_2}
\frac{\largetr_1'+\sigma_1x^{-1}}{\largetr_2'+\sigma_2x^{-1}}=
\exp(\largetr_1-\largetr_2)x^{\sigma_1-\sigma_2}(1+\smalltr)\gg 1$$

(2) We now let $\anytr$ be arbitrary with the property
$\mag(\anytr)\ne 1$ and  use  Proposition~\ref{Pstd} to write

$$\anytr=c\anymo_0+c\sum_{\bk\ge\bk_0}c_\bk\anymo_0\smallmo^\bk
=c\anymo+c\sum_{\bk\ge\bk_0}c_\bk\anymo_\bk$$

\z where $\anymo_\bk\ll \anymo_0$ and thus, by step (1) we have

$$\anytr'
=c\anymo'+c\sum_{\bk\ge\bk_0}c_\bk\anymo_\bk'=c\anymo'(1+\smalltr)$$

\z The rest of the proof is immediate.\end{proof}

\begin{Corollary}
  We have ${\mathcal{D}}\anytr=0\iff\anytr=\mathrm{Const.}$
\end{Corollary}
\begin{proof}
  We have to show that if $\anytr=\largetr+\smalltr\ne 0$ then
  $\anytr'\ne 0$. If $\largetr\ne 0$ then (for instance) $\largetr+\smalltr\gg
  x^{-1}=\smalltr$ and then $\largetr'+\smalltr'\gg x^{-2}\ne 0$. If
  instead $\largetr=0$ then $(1/\anytr)=\largetr_1+\smalltr_1+c$ and
we see that $(\largetr_1+\smalltr_1)'=0$ which, by the above, implies
$\largetr_1=0$ which gives $1/\smalltr=\smalltr_1$, a contradiction.
\end{proof}

\begin{Proposition}
  \label{P dom} Assume $\anytr=\largetr$ or $\anytr=\smalltr$. Then:

\z (i) If 
 $l(\mag(\anytr))\ge 1$ then
$l(\mag(\anytr^{-1}\anytr'))<l(\mag(\anytr))$.

\z (ii) 
$\dom(\anytr')= \dom(\anytr)'(1+\smalltr)$.

\end{Proposition}

\begin{proof}
  Straightforward induction.
\end{proof}

\subsubsection{Transseries with complex coefficients}\label{complextr}
Complex transseries $\mathcal{T}_{\CC}$ are constructed in a similar
way as real transseries, replacing everywhere $\largetr_1>\largetr_2$
by $\Re\largetr_1>\Re\largetr_2$. Thus there is only one order
relation in $\mathcal{T}_{\CC}$, $\gg$. Difficulties arise when
exponentiating transseries whose dominant term is imaginary.
Operations with complex transseries are then limited.  We will only
use complex transseries in contexts that will prevent these
difficulties.

\subsubsection{Differential systems in $\mathcal{T}$}

The theory of differential equations in $\mathcal{T}$ is similar to the
corresponding theory for functions.

{\em Example}.  The general solution of the differential equation 

\begin{equation}
  \label{eqEi1}
  f'+f=1/x
\end{equation}

\z in $\mathcal{T}$ (for $x\rightarrow +\infty$)
is $\anytr(x;C)=\sum_{k=0}^{\infty}k!x^{-k}+Ce^{-x}=\anytr(x;0)+Ce^{-x}$.

Indeed, the fact that $\anytr(x;C)$ is a solution
follows immediately from the definition of the operations
in $\mathcal{T}$. To show uniqueness, assume $\anytr_{1}$
satisfies  (\ref{eqEi1}). Then $\anytr_2=\anytr_{1}-\anytr(x;0)$
is a solution of ${\mathcal{D}}\anytr+\anytr=0$. Then 
$\anytr_2=e^{x}\anytr$ satisfies ${\mathcal{D}}\anytr_2=0$
i.e., $\anytr_2=\mathrm{Const.}$

The particular solution $\anytr(x;0)$ is the unique solution of the
equation $f=1/x-{\mathcal{D}} f$ which is manifestly
contractive\index{asymptotically contractive} in the space of level
zero transseries \index{transseries} (cf. \S~\ref{SubContr}). However
this same equation is not contractive\index{asymptotically
  contractive} for transseries \index{transseries} of positive level,
(because e.g. ${\mathcal{D}} e^x=e^x$); this could also have been
anticipated noting that the solution is not unique.

\subsubsection{Restricted composition}\label{Scompo}
The right composition $\anytr_1\circ\anytr_2$ is defined on
$\spacetr$, if $\mag(\anytr_2)\gg 1$ and $\dom(\anytr_2)>0$. The
definition is inductive.  

We first define the power and the exponential of a transseries \index{transseries}. Assume
powers and exponentials have been defined for all transseries \index{transseries} of level
$\le n-1$. Let
$\anytr=c\,\mag(\anytr)(1+\smalltr)\in\tilde{\mathcal{T}}^{[n]}$ be any
transseries \index{transseries} such that $ c>0$, cf.  Proposition~\ref{Pstd}.  By the
definition of $\mag(\cdot)$ and (\ref{defanytrn}), $\mag(\anytr)$ is a
transmonomial, $\mag(\anytr)=\largemo^{[n-1]}\exp(\largetr^{[n-1]})$.
We let

\begin{multline}
  \label{Edefpower}
  \anytr^\sigma=c^{\sigma}\left(\largemo^{[n-1]}\right)^{\sigma}
  \exp(\sigma\largetr^{[n-1]}) (1+\smalltr)^{\sigma}\\=
  c^{\sigma}\left(\largemo^{[n-1]}\right)^{\sigma}
  \exp(\sigma\largetr^{[n-1]}) (1+\smalltr)^{\sigma}\\=
  c^{\sigma}\left(\largemo^{[n-1]}\right)^{\sigma}
  \exp(\sigma\largetr^{[n-1]})
  \sum_{k=0}^{\infty}\binom{n}{\sigma}\smalltr^{n}
\end{multline}
\z  where   $\binom{n}{\sigma}$    are    the  generalized    binomial
coefficients,  the    infinite     sum   is      well   defined,    by
Proposition~\ref{infisum}, and thus  $\anytr^{\sigma}$ is well defined
as well.  Then, if $\boldsymbol\sigma\in(\RR^+)^M$
and if $\anytr^{[0]}=\sum_{\bk\ge\bk_0}c_\bk x^{-\boldsymbol\sigma\cdot\bk}$
is a level zero transseries \index{transseries}, we write

$$\anytr^{[0]}\circ \anytr=
\sum_{\bk\ge\bk_0}c_\bk (\anytr^{-1})^{\boldsymbol\sigma\cdot\bk}$$

\z which is well defined by Proposition~\ref{infisum} (ii) and
Proposition~\ref{Pstd}. We note that, under our assumptions for
$\anytr$, $\anytr^{[0]}\circ\anytr>0$ is positive iff $\anytr^{[0]}>0$

Similarly, we write cf. Definition~\ref{smalllarge}

\begin{multline}
  \exp(\anytr)=e^{L(\anytr)+C(\anytr) +s(\anytr)}\\=
e^{C(\anytr)}e^{L(\anytr)}\sum_{k=0}^{\infty}\frac{s(\anytr)^n}{n!}
=C'\anymo^{[n+1]}\anytr^{[n]}
\end{multline}

\z well defined by the definition of a transmonomial,
Proposition~\ref{infisum} (ii) and Proposition~\ref{Pstd}. Now the
definition of general composition is straightforward induction. We
assume that composition is defined at all $\le n-1$ levels, and that
in addition $\anytr^{[n-1]}\circ\anytr>0$ if $\anytr^{[n-1]}>0$. Then

\begin{multline}
  \anytr^{[n]}\circ\anytr=\sum_{\bk\ge\bk_0}c_{\bk}(\smallmo^{[n]}
  \circ\anytr)^{\bk}
  \\=\sum_{\bk\ge\bk_0}c_{\bk}(\anymo^{[n-1]}\circ\anytr)^{\bk}
  \exp(-\largetr^{[n-1]}\circ\anytr)\\
=\sum_{\bk\ge\bk_0}c_{\bk}(\anymo^{[n-1]}\circ\anytr)^{\bk}
  \left[\exp(-L(\largetr^{[n-1]}\circ\anytr)C\exp(-s(\largetr^{[n-1]}\circ\anytr)\right]\\
=\sum_{\bk\ge\bk_0}c'_{\bk}(\anymo^{[n-1]}\circ\anytr)^{\bk}
  \smallmo^{[n]}\exp(-s(\largetr^{[n-1]}\circ\anytr)
\end{multline}

\z and the last sum exists by the induction hypothesis and
Proposition~\ref{genclos}.

\subsection{The space $\mathcal{T}$ of general transseries}

We define 

\begin{align}
  L_n(x)=\mathop{\underbrace{\log\log ... \log(x)}}_{n\ times}\\
E_n(x)=\mathop{\underbrace{\exp\exp... \exp(x)}}_{n\ times}\\
\end{align}
\z with the convention $E_0(x)=L_0(x)=x$.

\z We write $\exp(\ln x)=x$ and then any log-free transseries \index{transseries} can be
written as $\anytr(x)=\anytr\circ E_n(L_n(x))$. This defines right
composition with $L_n$ in this trivial case, as $\anytr_1\circ L_n(x))=(\anytr\circ
E_n)\circ L_n(x):=\anytr(x)$. 

More generally, we define $\mathcal{T}$, the space of general transseries \index{transseries}, 
as a set of formal compositions

$$\mathcal{T}=\{\anytr\circ  L_n:\anytr\in\spacetr\}$$

\z with the algebraic operations (symbolized below by $*$) inherited
from $\spacetr$ by

\begin{align}
  \label{defop}
(\anytr_1\circ  L_n)*(\anytr_2\circ  L_{n+k})=
\left[(\anytr_1\circ E_k) *\anytr_2\right]\circ L_{n+k}
\end{align}

\z and using (\ref{defop}), differentiation is defined by 
$${\mathcal{D}}( \anytr\circ
L_n)=\left[(\prod_{k=0}^{n-1}L_k)^{-1}\right]({\mathcal{D}}\anytr)\circ L_n$$

\begin{Proposition}
  \label{Pbigfield}
$\mathcal{T}$ is an ordered differential field, closed under restricted composition.
\end{Proposition}

\begin{proof}
  The proof is straightforward, by substitution from the results in
  \S~\ref{Slogfree}.
\end{proof}

We will denote generically the elements of $\mathcal{T}$
with the same symbols that we used for $\spacetr$.
\begin{Proposition}
\label{Pintegr}
$\mathcal{T}$ is closed under  integration.

\end{Proposition}
\begin{proof}
 
The idea behind the construction of ${\mathcal{D}}^{-1}$ is the following: we
first find an  invertible operator $J$ which is {\em to leading order}
$\mathcal{D}^{-1}$; then the equation for the correction will be
contractive\index{asymptotically contractive}.
Let $\anytr=\sum_{\bk\ge\bk_0}\smallmo^\bk\circ L_n$. To unify 
the treatment, it is convenient to use the identity
$$\int_x\anytr(s)ds=\int_{L_{n+2}(x)}\left(\anytr\circ
E_{n+2}\right)(t)\prod_{j\le
  n+1}E_{j}(t)\mathrm{d}t\,=\int_{L_{n+2}(x)}\anytr_1(t)\mathrm{d}t\,$$

\z where the last integrand, $\anytr_1(t)$ is a log-free transseries \index{transseries} and
moreover 

$$\anytr_1(t)=\sum_{\bk\ge\bk_0}c_{\bk}\smallmo_1^{k_1}\cdots\smallmo_M^{k_M}=
\sum_{\bk\ge\bk_0}c_{\bk}e^{-k_1\largetr_1-...-k_M\largetr_M} $$

\z The case $\bk=0$ is trivial and it thus suffices to find
$\mathcal{D}^{-1} e^{-\largetr}$, where $n=l(\largetr)\ge 1$.  Then  $\largetr\gg x^m$ for any $m$ and thus also $\mathcal{D}\largetr\gg
x^m$ for all $m$. Therefore, since
$\mathcal{D}e^{-\largetr}=-(\mathcal{D}\largetr) e^{-\largetr}$ we expect that
$\dom(\mathcal{D}^{-1}e^{-\largetr})=-(\mathcal{D}\largetr)^{-1}e^{-\largetr}$
and we look for a $\Delta$ such that
\begin{equation}
  \label{defDelta}
\mathcal{D}^{-1}e^{-\largetr}=-\frac{e^{-\largetr}}{\mathcal{D}\largetr}(1+\Delta)  
\end{equation}

Then $\Delta$ satisfies the equation

\begin{equation}
  \label{eq:intContrac}
  \Delta=-\frac{\mathcal{D}^2\largetr}{(\mathcal{D}\largetr)^2}-
\frac{\mathcal{D}^2\largetr}{(\mathcal{D}\largetr)^2}\Delta+(\mathcal{D}\largetr)^{-1}\mathcal{D}\Delta
\end{equation}

\z By Propositions~\ref{Rform0},
Corollary~\ref{Cfinigen} and Proposition~\ref{P dom}, (\ref{eq:intContrac})
is contractive\index{asymptotically contractive} in $\mathcal{T}^{[n]}$. The proof now follows from  Proposition~\ref{genclos}.

\end{proof}

\z In the following we also use the notation ${\mathcal{D}}\anytr=\anytr'$
and we write $\mathcal{P}$ for the antiderivative $\mathcal{D}^{-1}$
constructed above.

\begin{Proposition}
  \label{noconst}
$\mathcal{P}$ is an antiderivative without constant terms, i.e,

$$\mathcal{P}\anytr=\largetr+\smalltr$$
\end{Proposition}

\begin{proof}
  This follows from Proposition~\ref{Rform0}, together with the fact that 
$\Delta,\largetr$ and $\mathcal{D}$ in (\ref{defDelta}) belong to $\mathcal{T}^{[n]}$.
\end{proof}
\begin{Proposition}
  \label{Ppropint}
We have 
\begin{align}
  \label{Epropint}
\mathcal{P}(\anytr_1+\anytr_2)=\mathcal{P}\anytr_1+\mathcal{P}\anytr_2\nonumber\\
(\mathcal{P}\anytr)'=\anytr;\ \mathcal{P}\anytr'=\anytr(0)\nonumber\\
\mathcal{P}(\anytr_1\anytr_2')=\anytr_1\anytr_2-\mathcal{P}(\anytr_1'\anytr_2)\nonumber\\
\anytr_1\gg\anytr_2\implies \mathcal{P}\anytr_1\gg \mathcal{P}\anytr_2\nonumber\\
\anytr>0\implies\mathcal{P}\anytr>0
\end{align}

\z where

 $$\anytr=\sum_{\bk\ge\bk_0}c_\bk\smallmo^{\bk}\implies\anytr(0)=\sum_{\bk\ge\bk_0;\bk\not =
  0}c_\bk\smallmo^{\bk}$$

\end{Proposition}
\begin{proof}
  All the properties are straightforward; preservation of inequalities
uses Proposition~\ref{Porrder}.
\end{proof}

\begin{Remark}
  \label{Rexpcontrac}
Let $\smalltr_0\in\mathcal{T}$. The operators defined by

\begin{align}
  \label{EJ1}
J_1(\anytr)=\mathcal{P}(e^{-x}(\mathrm{Const.}+\smalltr_0)\anytr(x))\\
J_2(\anytr)=e^{\pm x}x^{\sigma}\mathcal{P}(x^{-2}x^{-\sigma}e^{\mp
  x}(\mathrm{Const.}+\smalltr_0)\anytr(x))\label{EJ2}
\end{align}

\z are contractive\index{asymptotically contractive} on $\mathcal{T}$.
\end{Remark}

\begin{proof} For (\ref{EJ1}) it is   enough to show contractivity\index{asymptotically contractive}
  of $\mathcal{P}(e^{-x}\cdot)$. This is a straightforward calculation
  similar to the proof of Proposition~\ref{Pintegr}.  We have for some
  $n$ $\anytr(x)=\sum_{\bk\ge\bk_0}\smallmo^{\bk}(L_n(x))$ where
  $\smallmo_j\in\stackrel{\approx}{\mathcal{T}}$.
\begin{multline}
  \label{Eintg}
\mathcal{P}e^{-x}(\anymo\circ L_n)=\mathcal{P}\left(e^{-E_{n+2}}\prod_{1\le
    j\le
    n+2}E_j\exp(\largetr\circ E_2)\right)\circ L_{n+2}\\=
\left[\frac{e^{-E_{n+2}}\prod_{1\le
    j\le
    n+2}E_j\exp(\largetr\circ E_2)}{-E_{n+2}'+\sum_{0\le
    j\le
    n+1}E'_j+\largetr'\circ E_2'}\left(1+\smalltr\right)\right]\circ L_{n+2}\\ \ll
\prod_{1\le
    j\le
    n+2}E_j\exp(\largetr\circ E_2)
\end{multline}
  
\z The proof of (ii) is similar.

\end{proof}
\section{Equations
  in $\mathcal{T}$: examples}
\begin{Remark}
  The general contractivity principle stated in Theorem~\ref{PFp},
  which we have used in proving closure of transseries with respect to
  a number of operations can be used to show closure under more
  general equations. Our main focus is on differential systems.
\end{Remark}
\subsubsection{Nonlinear ODEs in $\mathcal{T}$}
\label{deg1}
We start with an example, a first order equation:

\begin{equation}
  \label{Eqord1}
f'=J_1(f)=F_0(x^{-1})-f-\frac{\beta}{x}f-g(x^{-1},f)  
\end{equation}

\z where 

\begin{align}
  \label{Edefterms}
F_0(x^{-1})=\sum_{k\ge 2}\frac{F_{0k}}{x^k}\nonumber\\
g(x^{-1},f)=\sum_{k\ge 0;\ l\ge 1}g_{kl}x^{-k}f^l
\end{align}

\z where the sums are assumed to converge and $g_{01}=g_{11}=0$.

\z We see that $J_1$ is well defined if $f=\smalltr\in\mathcal{T}$
(cf. Proposition~\ref{genclos}), and it is under this assumption that
we study $J_1$.\footnote{If there are infinitely many nonzero terms in
  the sum in (\ref{Edefterms}), $J_1$ is not in $\mathcal{T}$ if $f\gg
  1$ (since, in this case, $\mag(f^n)$ is unbounded).}

(1). Solutions of (\ref{Edefterms}) in $\tilde{\mathcal{T}}^{[0]}$.
The equation 

\begin{eqnarray}
  \label{Elev0}
f=J_2(f)=-{\mathcal{D}} f+F_0(x^{-1})-\frac{\beta}{x}f-g(x^{-1},f)  
\end{eqnarray}

\z is contractive\index{asymptotically contractive} in $\tilde{\mathcal{T}}^{[0]}$ (this follows
immediately from \S\ref{SubContr}). Thus there exists in
$\tilde{\mathcal{T}}^{[0]}$ a unique solution $\tilde{f}_0$. Since 
(\ref{Elev0}) is also contractive\index{asymptotically contractive} in the subspace of
$\tilde{\mathcal{T}}^{[0]}$
of series of the form $\sum_{k=2}^{\infty}\frac{c_k}{x^k}$ we have

\begin{eqnarray}
  \label{Eformy0}
  \tilde{f}_0=\sum_{k=2}^{\infty}\frac{c_k}{x^k}
\end{eqnarray}

{\bf Note.}  The iteration $f_{n+1}=J_1 f_n$, $f_1=x^{-1}$ is convergent
in $\mathcal{T}$ and, if $f_i=\sum_{k=2}^{i}c_k^{[i]}x^{-k}$ then
$c_k^{[i]}=c_k$ for $k\le i$, and this is a very
convenient way to calculate the coefficients $c_i$.

(2)  Let now $\delta=f-\tilde{f}_0$. Then

\begin{align}
  \label{Eeqdelta}
\delta'=-\delta-\frac{\beta}{x}\delta-g(x^{-1},\tilde{f}_0+\delta)
+g(x^{-1},\tilde{f}_0+\delta)\nonumber\\
=-\delta-\frac{\beta}{x}\delta+\sum_{k\ge 0;\ l\ge 1}c_{kl}x^{-k}\delta^l
\end{align}
\z with 
\begin{align}
  \label{Ecoeffdelta}
c_{01}=c_{11}=0
\end{align}

\z or

\begin{align}
  \label{Eeqdelta1}
\frac{\delta'}{\delta}=-1-\frac{\beta}{x}-\sum_{k\ge
  2}\frac{c_{k1}}{x^k}+
\sum_{k\ge 0;\ l\ge 1}c_{k;l+1}x^{-k}\delta^l
\end{align}

\z Since by assumption $\delta\ll 1$ we have
$$\ln\delta=C_0-x+\beta\ln x+\sum_{k\ge
  1}\frac{c_{k+1;1}}{kx^k} + x\smalltr (x)$$

\z and thus $\delta\ll \exp(-cx)$ for any $c<1$ so that 

$$\ln\delta=C_0-x+\beta\ln x+\sum_{k\ge
  1}\frac{c_{k+1;1}}{kx^k}+\exp(-cx)\smalltr (x)$$

\z whence, by composition with $\exp$ we get

$$\delta=C_1x^{\beta}e^{-x}\sum_{k\ge
  1}\frac{d_{k+1;1}}{kx^k}+\exp(-cx)\smalltr (x)$$

\z Equation (\ref{Eeqdelta1}) implies 

\begin{align}\label{Eequint}
  \delta=Cx^{\beta}e^{-x}\tilde{y}_0\exp\left(\int\sum_{k\ge 0;\ l\ge
    1}c_{k;l+1}x^{-k}\delta^l\right);\ \ \left(\tilde{y}_0=\sum_{k\ge
  0}\frac{d_{k+1;1}}{kx^k}\right)
\end{align}

\z and (\ref{Eequint}) is contractive\index{asymptotically contractive} by Remark~\ref{Rlinearity} and
Remark~\ref{Rexpcontrac}.  In particular, for every $C$ there is a
unique $\delta(x;C)$ satisfying (\ref{Eequint}).

\begin{Proposition}
  \label{Rspecialform}
The general transseries solution of (\ref{Eqord1}) is $\tilde{f_0}+\delta$
where  
\begin{align}
 \label{Especialform} 
\delta=\sum_{k=1}^{\infty}C^kx^{\beta k}e^{-kx}\tilde{f}_k(x)
\end{align}

\z with $\tilde{f}_k\in\tilde{\mathcal{T}}^{[0]}$ and

$$\tilde{f}_k(x)=\sum_{j=0}^{\infty}\frac{f_{k;j}}{x^j}$$
\end{Proposition}
\begin{proof}
  This is a straightforward consequence of discussion of this section
  and of (\ref{Eequint}) .
\end{proof}

\subsubsection{Formal linearization}
\label{sec:Lineariz}
Let $z=Cx^{\beta}e^{-x}$. We have $$C(x,\delta)=x^{-\beta}e^{x}\sum_{k\ge
  1}\delta^k\tilde{g}_k(x)$$ A direct calculation
shows that $C'=C_x+C_\delta\delta'=0$. The transformation
$(x\mapsto x;y\mapsto C(x,y-f_0)) $ formally linearizes (\ref{Eqord1}).

\subsection{Multidimensional systems: transseries solutions at irregular
  singularities of rank one}\label{Sysnonlin} Consider the differential
system
            
\begin{eqnarray}
 \label{eqor1}
  \mathbf{y}'=\mathbf{f}(x^{-1},\mathbf{y})  \qquad \mathbf{y}\in\CC^n               \end{eqnarray}

   We look at solutions $\mathbf{y}$ such that $\mathbf{y}(x)\rightarrow
   0$ as $x\rightarrow\infty$ along some direction
   $d=\{x\in\CC:\arg(x)=\phi\}$. The
   following conditions are assumed

\z (a1) The function $\mathbf{f}$ is analytic
at $(0,0)$.

\z (a2) Nonresonance: the eigenvalues $\lambda_i$ of the linearization

\begin{eqnarray}
  \label{linearized}
  \hat{\Lambda}:=-\left(\frac{\partial f_i}{\partial
    y_j}(0,0)\right)_{i,j=1,2,\ldots n}
\end{eqnarray}

\z are linearly independent over $\ZZ$ (in particular nonzero) and such
that $\arg\lambda_i$ are different from each other (i.e., the Stokes
lines are distinct; we will require somewhat less restrictive
conditions, see \S~3.1). 

By relatively straightforward algebra, (\cite{Wasow} and also \S~3.2
where all this is exemplified in a two-dimensional case), the system
(\ref{eqor1}) can then be brought to the form

\begin{eqnarray}\label{eqor}
{\bf y}'=-\hat\Lambda {\bf y}+ \frac{1}{x}\hat A {\bf y}+{\bf
g}(x^{-1},{\bf y})
\end{eqnarray}

\z where $\hat{\Lambda}=\mbox{diag}\{\lambda_i\},\ 
\hat{A}=\mbox{diag}\{\alpha_i\}$ are constant matrices, $${\bf
  g}(x^{-1},{\bf y})= O(x^{-2},\mathbf{y}^2),\ 
(x\rightarrow\infty,\mathbf{y}\rightarrow 0)$$
\begin{Remark}\label{invert}
  (i) If ${\bf g}(x^{-1},{\bf y})\equiv 0$. In this case the system
  (\ref{eqor}) is linear and has the general transseries solution

$${\bf y}=\erm^{-x\hat{\Lambda}}\mathbf{C}
  x^{\hat{A}}$$ 

(ii) More generally, if ${\bf g}(x^{-1},{\bf y})={\bf G}(x^{-1})$ is a
transseries,
then the general solution of (\ref{eqor}) is

\begin{equation}
  \label{eq:inhom}
  {\bf y}=\erm^{-x\hat{\Lambda}}
  x^{\hat{A}}\mathbf{C}+\erm^{-x\hat{\Lambda}}
  x^{\hat{A}}\mathcal{P}\Big(\erm^{x\hat{\Lambda}} x^{-\hat{A}}{\bf
    g}\Big)
\end{equation}

\end{Remark}

\begin{proof}
  In both cases the system is diagonal and the result follows
  immediately from the case when $n=1$, i.e. from
  Proposition~\ref{Rspecialform}.
\end{proof}

The general solution of (\ref{eqor}) in $\mathcal{T}_{\CC}$  (cf.
\S\ref{complextr}) is an
$n_1\le n$ parameter transseries, as shown in the sequel.

\begin{Proposition}
  \label{GTS}
  
  Let $d$ be a ray in $\CC$.  The general solution of (\ref{eqor}) in
  $\mathcal{T}_\CC$ with the restriction $\mathbf{y}\ll 1$ is of the
  form

\begin{gather}
  \label{transsf} \tilde{\mathbf{y}}(x)=\sum_{\mathbf{k}\ge 0}
  \mathbf{C}^{\mathbf{k}}\erm^{-\boldsymbol{\lambda}\cdot\mathbf{k}x}
  x^{\boldsymbol{\alpha}\cdot\mathbf{k}}\tilde{\mathbf{s}}_{\mathbf{k}}(x)=
  \sum_{\mathbf{k}\ge 0}
  \mathbf{C}^{\mathbf{k}}\erm^{-\boldsymbol{\lambda}\cdot\mathbf{k}x}
  x^{\mathbf{m}_0\cdot\mathbf{k}}\tilde{\mathbf{y}}_{\mathbf{k}}(x)
\end{gather}

\end{Proposition}
\z where $C_i= 0$ for all $i$ so that $\erm^{-\lambda_i
  x}\not\rightarrow 0$ as $x\rightarrow\infty$ in $d$.

\begin{proof}
If $\mathbf{y}$ is a solution of (\ref{eqor}) then we have, by
Remark~\ref{invert}

\begin{equation}
  \label{eq:invert}
  \mathbf{y}=\erm^{-x\hat{\Lambda}}
  x^{\hat{A}}\mathbf{C}+\erm^{-x\hat{\Lambda}}
  x^{\hat{A}}\mathcal{P}\Big(\erm^{x\hat{\Lambda}} x^{-\hat{A}}\mathbf
    {g}(x^{-1},\mathbf{y})\Big)
\end{equation}

\z for some $\mathbf{C}$. Since $\mathbf{y}\ll 1$ we have
$\mathbf{g}(x^{-1},\mathbf{y})\ll 1$ and thus $$\mathcal{P}\Big(\erm^{x\hat{\Lambda}} x^{-\hat{A}}\mathbf
    {g}(x^{-1},\mathbf{y})\Big)\ll \erm^{x\hat{\Lambda}} x^{-\hat{A}}$$ 
Again since $\mathbf{y}\ll 1$, we then have $C_i=0$ for all 
$i$ for which $e^{-\lambda_i x}\not\ll 1$.

\z {\bf Note}.  With the condition $\mathbf{y}\ll 1$,
eq. (\ref{eq:invert}) has a unique solution. 

Indeed, the difference of
two solutions $\mathbf{y}_1 -\mathbf{y}_2$ satisfies the equation
\begin{equation}
  \label{eq:invert2}
 \mathbf{y}_1
-\mathbf{y}_2=\erm^{-x\hat{\Lambda}}
  x^{\hat{A}}\mathcal{P}\Big(\erm^{x\hat{\Lambda}} x^{-\hat{A}}\big[\mathbf
    {g}(x^{-1},\mathbf{y}_1)-\mathbf
    {g}(x^{-1},\mathbf{y}_2)\big]\Big)
\end{equation}

\z Since $\mathbf {g}(x^{-1},\mathbf{y})=O(x^{-2},\mathbf{y}^2)$ we
    have $$\mathbf {g}(x^{-1},\mathbf{y}_1)-\mathbf
    {g}(x^{-1},\mathbf{y}_2)=O(x^{-2}\boldsymbol{\delta},|\mathbf{y}||\boldsymbol{
    \delta}|)$$ which by Proposition~\ref{Ppropint} implies
    $\boldsymbol{\delta}=o(\boldsymbol{\delta})$, i.e., $\delta=0$.

Using Remark~\ref{Rexpcontrac} it is easy to check that
(\ref{eq:invert2}) is an asymptotically contractive\index{asymptotically contractive} equation, in the
space of $\mathbf{y}$ which are $\ll x^{-2}$ thus it has a solution
$\mathbf{y}^{[0]}$ with this property. Since the previous note shows the
solution of (\ref{eqor}) with $\mathbf{y}\ll 1$ is unique, we have
$\mathbf{y}=\mathbf{y}^{[0]}$. Formula (\ref{transsf}) is obtained
by straightforward iteration of (\ref{eq:invert2}).

\end{proof}



\section{Borel summation techniques} 
\subsubsection{Introduction}

In this section we discuss a Borel summation--induced isomorphism
between transseries and functions in the setting \S\ref{Sysnonlin}.
The goal is to show, on simple examples, how Borel summation is proved
and used.  For more general results we refer to \cite{Duke}.

\begin{Definition}
A Borel-summable series $\tilde y:=\sum_{k=K}^\infty y_kx^{-k}$,
$K\in\ZZ$ is a formal
power series with the following properties
\begin{enumerate}
\item the truncated Borel transform $Y=\mathcal{B}\tilde
y:=\sum_{k>0}\frac{y_{k}}{(k-1)!}t^{k-1}$ of $\tilde y$ has a
nonzero radius of convergence, \item $Y$ can be analytically
continued along $[0,+\infty)$ and \item the analytic continuation
$Y$ grows at most exponentially along $[0,+\infty)$ and is
therefore Laplace transformable along $[0,+\infty)$.
\end{enumerate}
The Borel sum $y$ of $\tilde y$ is then given by
\begin{equation}
  \label{eq:defLB}
  y = \lap\bor\tilde{y}:=\sum_{k=K}^{0}y_{k}x^{-k} +
  \mathcal{L}Y,
\end{equation}
where the sum is understood to be zero if $K>0$ and $\mathcal{L}$
denotes the usual Laplace transform.
\end{Definition}
{\em Example} The formal solution for large $x$ of the equation
$f'-f=x^{-1}$ is Borel summable. Indeed
$$\lap\bor\sum_{k=0}^{\infty}\frac{-k!}{(-x)^{k+1}}\\=\lap\frac{1}{1+p}=\int_0^{\infty}\frac{e^{-px}}{1+p}dp$$
and it can be checked that the Borel sum is a solution of the given
equation; see also (i) in the note below. 
\begin{Note}{\rm 
    (i) It can be shown that Borel summation is an extended
    isomorphism: in particular it commutes with algebraic operations,
    including multiplication and complex conjugation, and with
    differentiation \cite{Balser}. This can be expected from the fact
    that, formally, Borel summation is the composition of Laplace
    transform with its inverse, and the identity obviously commutes
    with the operations mentioned above.
  
  (ii) Borel summation has to be generalized in a number of ways in
  order to be used for solving more general equations. Indeed, an
  equation as simple as $f'+f=x^{-1}$ has a formal solution which is
  not Borel summable: $\bor \sum_{k=0}^{\infty}k! x^{-k-1}=(1-p)^{-1}$
  is not Laplace transformable. If the Laplace transform integral is
  taken along a ray above or below $\RR^+$, then $\lap\bor (1-p)^{-1}$
  is a solution of the given equation. None of these path choices
  however yields a real valued function, whereas the formal series has
  real coefficients. The resulting summation operator would not be an
  isomorphism since commutation with complex conjugation fails. In
  this simple example the half sum of the upper and lower integrals is
  real valued and solves the equation, but this ad hoc procedure
would not work for nonlinear equations. More complicated averages 
need to be introduced, see \S\ref{Ave}. }
\end{Note}
 \subsection{Borel summation of transseries: a first order example}\label{Ord1}
 
 Consider the (first order) differential equation (\ref{Eqord1}) with,
 say, $F_0=x^{-2}$ and $g=a f^2+b f^3$:

 \begin{equation}
   \label{eq:eq*}
   f'+(1-\beta x^{-1})f=x^{-2}+af^2+bf^3
 \end{equation}

\z  and $a,b$ some constants. We first look at solutions, both
 formal and actual, which go to zero as $x\rightarrow\infty$ in some
 direction in $\CC$.

We have seen in \S~\ref{deg1} that the solution $\tilde{f}_0$ as a
 formal power series (i.e., in $\mathcal{T}^{[0]}$), is unique and
 $\tilde{f}_0=\sum_{k=2}^{\infty}c_kx^{-k}$. By Remark~\ref{Rspecialform}
the general transseries solution is

\begin{equation}
  \label{eq:gtrns}
  \tilde{f}=\sum_{k=0}^{\infty}\xi^k\tilde{f}_k(x)\ \ \ \ (\xi:=Ce^{-x}x^\beta)
\end{equation}

\z where $\tilde{f}_k$ are integer power series.  It is important to
note that in (\ref{eq:gtrns}) only the constant $C$ depends on the
solution. We will see that the $\tilde{f}_k$ are simultaneously Borel
summable, to $f_k=\mathcal{LB}\tilde{f}_k$, and that the sum

\begin{equation}
  \label{eq:gtrns,sum}
 f=\sum_{k=0}^{\infty}\xi^k{f}_k(x)
\end{equation}

\z is {\em convergent} and provides the general solution of
(\ref{eq:eq*}) with the property $f\rightarrow 0$ as
$x\rightarrow\infty$ in some direction in which (\ref{eq:gtrns}) is a
valid complex transseries.

Assuming for the moment we proved that (\ref{eq:gtrns,sum}) indeed
provides a solution of (\ref{eq:eq*}) for any $C$ it is not difficult
to show that there are no further solutions:

\begin{Lemma}
  \label{uniquen}
If $f_1$ is any solution of (\ref{eq:eq*}) with the stated condition 
in some direction at
infinity, 
then $f_1-f_0=Ce^{-x}x^\beta(1+o(1))$ as $x\rightarrow\infty$.
If $C=0$ then $f_1=f_2$.

\end{Lemma}

\begin{proof}
 Writing the equation for $\delta=f_1-f_0$, multiplying
with the integrating factor  of the ``dominant''
part of the equation, $Ce^{t}t^{-\beta}$, and integrating
we get

 \begin{equation}
   \label{eq:intequ}
   \delta = Ce^{-x}x^\beta+ e^{-x}x^\beta\int_{a}^x
e^t t^{-\beta}\left[\left(2af_0 +3bf_0^2\right)\delta+(a+3b f_0)\delta^2+\delta^3\right]dt
 \end{equation}
 
 \z which is contractive in the sup norm on $(x_0,\infty)$, for $|x_0|$
 large enough, in a ball of radius $\epsilon>0$ small enough.
The solution of (\ref{eq:intequ}) is thus unique and it is easy
to see that $\delta =C_1e^{-x}x^\beta(1+o(1))$ for large $x$
($C_1$ is not in general equal to $C$).
\end{proof}

\z 
\begin{Corollary}
  Formula (\ref{eq:gtrns,sum}) provides the most general solution
of (\ref{eq:eq*}) with the property $f\rightarrow 0$ in
some direction in $\CC$.

\end{Corollary}

 A convenient way to
 generate $\tilde{f}_0$ is the iteration(\ref{Elev0}), which we start
 with $\tilde{f}_0^{[0]}=0$.  Denoting
 $\delta^{[k]}=\tilde{f}_0^{[k+1]}-\tilde{f}_0^{[k]}$ we have

$$\delta^{[k]}=(-\mathcal{D}-\beta x^{-1}+O(x^{-2}))\delta^{[k-1]}$$

\z whence $\delta^{[k]}\sim const. \Gamma(k-\beta)x^{-k}$ and thus
Borel summation appears natural.  Since $\tilde{f}_0$ is defined
through a differential equation it is natural to Borel transform the
equation itself. The result is

\begin{equation}
  \label{eq:{Eqord1,B}}
  -pF+F=p-\beta F*1+a F^{*2}+b F^{*3}
\end{equation}

\z where convolution is defined by
\begin{equation}
  \label{eq:defConvo}
(  f*g)(p)=\int_0^pf(s)g(p-s)ds
\end{equation}
and we write
$$F^{*k}=\underbrace{F*F*\cdots*F}_{\mbox{k times}}$$

\z  A solution 
$F$ that is Laplace transformable along any ray {\em other than
$\RR^+$} is obtained by noting that the equation  
\begin{equation}
  \label{eq:{Eqord1,B1}}
 F=(1-p)^{-1}\left(p-\beta F*1+a F^{*2}+b F^{*3})\right)=\mathcal{N}(F)
\end{equation}

\z is contractive in an appropriate norm.

\begin{Proposition}
  \label{C*a}
  The space $L^1_\nu$ of  functions along a ray 
$d=\{p:\arg(p)=a$, such that $\|f\|_\nu<\infty$ with
  the norm $\|f\|_\nu=\int_{t\in d} |f(t)|e^{-\nu |t|}d|t|$
  is a Banach algebra with respect to convolution.
\end{Proposition}

\begin{proof}
  All the properties are verified in a straightforward manner. 
In particular, we have $\|f*g\|_\nu=\|f\|_\nu \|g\|_\nu$. 
\end{proof}

\z \z Therefore, with $\rho=\|F\|_\nu$ and $d_1=\mbox{dist}(1,d)\ne 0$
we have

\begin{multline}
 \left\| \mathcal{N}(F)\right\|_\nu =\left\|(1-p)^{-1}\left(p-\beta
F*1+a F^{*2}+b F^{*3}\right)\right\|_\nu
\\\le \frac{1}{d_1}\left(\|p\|_\nu+|\beta| \rho
  \|1\|_\nu+ |a|\rho^2+ |b|\rho^3\right)\\=
\frac{1}{d_1}\left(\frac{1}{\nu^2}+\frac{|\beta|\rho}{\nu}
+ |a|\rho^2+ |b|\rho^3\right)\rightarrow 0 \ \ \mbox{as}\ \ \nu\rightarrow\infty
\end{multline}

\z and clearly, if $\nu$ and $1/\epsilon$ are large enough, the image
$\mathcal{N}B_\epsilon$
of the ball $B_\epsilon=\{F:\|F\|_\nu<\epsilon\}$ is contained 
in $B_\epsilon$. Similarly, it can be seen that $\mathcal{N}$ is contractive
in $B_\epsilon$. Thus the following conclusion.

\begin{Proposition}
  \label{L1loc}
There is a unique solution of
(\ref{eq:{Eqord1,B}}) in $\cup_{\nu\ge\nu_0}L^1_{\nu}$.

\end{Proposition}

This however does not yet imply Borel summability of the series
$\tilde{f}_0$; to show this we must prove appropriate analyticity 
properties for $F$ and this can be done in essentially the same manner. 

\begin{Proposition}
 \label{C*a1}
 The space of analytic functions in a region of the form

$$S_M=\{p:\arg(p)\in (a_1,a_2))\not\ni 0\ \mbox{and }|p|<M\}\cup\{p:|p|<1-\epsilon\}$$

\z vanishing at $p=0$, continuous in $\overline{S}$, and such that
$\|f\|_\nu<\infty$ with the norm
$\|f\|_{\nu;\infty}=M^{-1}\sup_{\overline{S}}|f(p)|e^{-\nu|p|}$ is a
Banach algebra with respect to convolution. In addition, if $g(\cdot
e^{i\phi})\in L^1_{\nu}(\RR^+)$ for any $\phi\in (a_1,a_2)$ and $g$ is
analytic in $S$ then we have

 \begin{align}
  \|f*g\|_{\nu;\infty} \le \|f\|_{\nu;\infty}\|g\|_{\nu;1}\\
 \end{align}

\z The function $F_0$ is analytic in $\mathcal{S}_M$,
 Laplace transformable in any direction other
than $\RR^+$ and $y_0=\lap \{F_0\}$ is a solution of (\ref{eq:eq*}).
\end{Proposition}

\begin{proof}
  For large enough $\nu$ and for any $M$ there is a unique analytic
  solution in $S_M$, $F_0$, and $F_0$ is thus independent of $M$. Since
  we have $|F_0(p)|\le |p|e^{\nu |p|}$ for $p\in S=\cup_{M>0} S_M$ it
  follows that $F_0=F$.  Using Proposition~\ref{L1loc} the proof is
  complete.
\end{proof}

\subsubsection{Summability of the transseries $\tilde{f}$}

A straightforward calculation shows that the series
$\tilde{f}_k$ in (\ref{eq:gtrns}) satisfy the system of equations

\begin{eqnarray}
  \label{eq:syste}
  \tilde{f}_0'+(1-\beta x^{-1})\tilde{f}_0 & = & x^{-2}+a\tilde{f}_0^2+b\tilde{f}_0^3\nonumber\\
  \tilde{f}_1'+(2a \tilde{f}_0 + 3 b \tilde{f}_0^2)\tilde{f}_1&=&0\\
 \tilde{f}_k'+\Big((1-k)(1-\beta x^{-1})+2a \tilde{f}_0 + 3 b
  \tilde{f}_0^2\Big)\tilde{f}_k
  &=&a{\sum}^\dagger\tilde{f}_{k_1}\tilde{f}_{k_2}
  +b{\sum}^{\ddagger}\tilde{f}_{k_1}\tilde{f}_{k_2}\tilde{f}_{k_3}\nonumber
\end{eqnarray}
 where in the sum $\sum^\dagger$ the indices satisfy $k_1+k_2=k;\ k_1>0;
k_2>0$ while in $\sum^{\ddagger}$ the condition is $ k_1+k_2+k_3=k ;\ 
k_i>0$. We have already seen that $\tilde{f}_0$ is summable. The
equation for $\tilde{f}_1$ is special, and we treat it separately. To
ensure Borel transformability, since $\tilde{f}_1=O(1)$, it is
convenient to take $\tilde{f}_1= x^2\tilde{y}_1$, which gives

$$\tilde{y}_1'+2x^{-1}\tilde{y}_1+(2a \tilde{f}_0 + 3 b
\tilde{f}_0^2)\tilde{y}_1=0$$

\z which, in Borel transform, with $Y_1= \bor\tilde{y}_1$,

\begin{equation}
  \label{eq:eqY10}
  -pY_1 +2\int_0^p Y_1(s)ds+(2a _0 F_0 + 3 b
F_0^{*2})*Y_1=0
\end{equation}

\z which implies, denoting $2a _0 F_0 + 3 b
F_0^{*2}=G$

\begin{equation}
  \label{eq:eqY1}
  -pY_1'+Y_1 =-G'*Y_1
\end{equation}

\z in which the leading behavior of $Y_1$ is expected to be
$Y_1=p+\ldots$ and thus the dominanat balance is between the terms on
the l.h.s. of (\ref{eq:eqY1}). In integral form we have, with $Y_1=pQ$,

\begin{equation}
  \label{eq:eqY11}
 Q=1+\int_0^p ds s^{-2}\int_0^suG'(u)Q(s-u)du=
1+\int_0^p\int_0^1   v G'(sv)Q(sv-v) dv ds
\end{equation}

\z It is easy to see that (\ref{eq:eqY1}) is contractive 
in a space of analytic functions for $|p|<\epsilon$
if $\epsilon$ is small. To find exponential bounds for large $p$
we first restrict to a ray $p=t e^{i\phi}$ with $\phi\ne 0$.

For $x=p e^{-i\phi}\in [0,\infty)$ it is useful to write $Q=Q_0+Q_1$
where $Q_0=0$ for $|p|>\epsilon$ and $Q_1=0$ for $|p|<\epsilon$. Noting
that $\int_0^s Q_1(t)G'(s-t)dt=0$ if $|s|<\epsilon$, the equation for
$Q_1$ takes the form, for $|p|>\epsilon$, 

$$Q_1=F(p)+\int_{\epsilon e^{i\phi}}^p 
\int_{\epsilon e^{i\phi}}^s s^{-2} G'(u)Q_1(s-u)\,\,\drm u\,\drm s $$

\z where 

$$F(p)=1+\int_0^p\int_0^{\epsilon e^{i\phi}}
s^{-2}G'(s-u)Q_0(u)\,\drm u\,\drm s
$$

\z Taking $Y(p)=Q_1(\epsilon e^{i\phi} +p)$ we get a convolution
equation
of the form

$$Y=F_1+Y*\{(s+\epsilon e^{i\phi})^{-2}\}*F_2$$

\z which is manifestly contractive in the norm $\|\cdot\|_\nu$
for large $\nu$. Since $Q_0$ is manifestly in $\mathcal{A}_\nu$
for any $\nu$, it follows that $Q\in\mathcal{A}_\nu$ as well.

For $k>1$ it is convenient to take $\tilde{f}_k=x^{2k}\tilde{y}_k$ 
and we have 

\begin{multline}\label{eqeqk}
  \tilde{y}_k'+\Big(1-k+(2k+(k-1)\beta) x^{-1}+2a \tilde{y}_0 + 3 b \tilde{y}_0^2\Big)\tilde{y}_k
=a{\sum}^\dagger
\tilde{y}_{k_1}\tilde{y}_{k_2}
\\+b{\sum}^{\ddagger}\tilde{y}_{k_1}\tilde{y}_{k_2}\tilde{y}_{k_3}
\end{multline}
(cf. (\ref{eq:syste})) which, after Borel transform becomes

\begin{multline}
  \label{eqeqYk}
  Y_k=(p+k-1)^{-1}\Bigg((2k+k\beta -\beta+G)*Y_k+
 a{\sum}^\dagger
 Y_{k_1}* Y_{k_2}
+b{\sum}^\ddagger
 Y_{k_1}* Y_{k_2}* Y_{k_3}\Bigg)\\=
\mathcal{J}Y_k+(p+k-1)^{-1}a{\sum}^\dagger
 Y_{k_1}* Y_{k_2}+b{\sum}^{\ddagger}
 Y_{k_1}* Y_{k_2}* Y_{k_3}
\end{multline}
Note now that the operator on the r.h.s. of 
(\ref{eqeqYk}) is contractive in the norms introduced, for large $\nu$. 
It is then clear that $Y_k$ are analytic and Laplace transformable
for large enough $\nu$. It remains to show they are
simultaneously Laplace transformable.

Note that for large $\nu$, denoting $\|Y_k\|=x_k$ we have, with
$\lambda_\nu=\max\{\|\mathcal{J}\|,\|Y_1\|_\nu\}$ arbitrarily small if
$\nu$ is large,
\begin{equation}
  \label{eqeqYk1}
  x_k\le \lambda_\nu x_k+a{\sum}^\dagger
x_{k_1}x_{k_2}+b{\sum}^\ddagger
 x_{k_1}x_{k_2}x_{k_3}
\end{equation}and the coefficients $x_k$ are majorized by the Taylor coefficients
of the analytic solution of the algebraic equation
$$\psi(z)=\lambda_\nu z+\lambda_\nu\psi(z)+a\psi(z)^2+b\psi(z)^3$$
Thus, for large enough $\nu$ we have that $x_k\le \rho_\nu^k$
($\rho_\nu=o(1)$ for large $\nu$) and thus, for $x$ large enough
$$\|\lap\{Y_k\}\|_{\infty}\le \rho_\nu^k$$
  (An inductive calculation
shows that $Y_k=O(p^{2k-1})$ for small $p$.) The sum
$$\sum_{k=0}^{\infty}(\xi x^2)^k \lap\{Y_k\}$$
is then uniformly convergent. It is easy to see that it gives
therefore a solution of (\ref{eq:eq*}). 
\subsection{Generalized Borel summation for rank one ODEs}

 We look at the differential system (\ref{eqor1}) under the same
assumptions and normalization as in \S~\ref{Sysnonlin}.

\z {\em Further normalizing transformations}.  For convenience, we
rescale $x$ and reorder the components of $\mathbf{y}$ so that

(n3) $\lambda_1=1$, and, with $\phi_i=\arg(\lambda_i)$, we have
$\phi_i<\phi_j$ if $i<j$. To simplify notations, we formulate some of
our results relative to $\lambda_1$; they can be easily adapted to any
other eigenvalue.

  To unify the treatment we
make, by taking $\mathbf{y}=\mathbf{y}_1 x^{-N}$ for some $N>0$,

(n4) $\Re(\beta_j)<0,\ j=1,2,\ldots,n$.

\z {\bf Note}: there is an asymmetry at this point: the opposite
inequality cannot be achieved, in general, as simply and without
violating analyticity at infinity. In some instances this
transformation is not convenient since it makes more difficult the
study of certain properties, see \cite{Invent}.

 Finally, through a transformation of the form
$\mathbf{y}\leftrightarrow\mathbf{y}-\sum_{k=1}^M\mathbf{a}_k x^{-k}$ we arrange that

(n5) $ \mathbf{f}_0=O(x^{-M-1})\mbox{ and }\mathbf{g}(x,\mathbf{y})=
O(\mathbf{y}^2,x^{-M-1}\mathbf{y}) $. We choose $M>1+\max_i\Re(-\beta_i)$.

{\em Formal solutions.} The transseries solutions of (\ref{eqor1}) 
were studied in \S~\ref{Sysnonlin}. More generally, there 
is an n--paramter family of formal exponential series
solutions of (\ref{eqor1}):
\begin{eqnarray}
  \label{eqformgen,n}
   \tilde{\mathbf{y}}_0+\sum_{\mathbf{k}\ge 0; |\mathbf{k}|>0}C_1^{k_1}\cdots C_n^{k_n}
\mathrm{e}^{-(\bfk\cdot\bflam) x}x^{\bfk\cdot\bfm}\tilde{\mathbf{y}}_{\bfk}
\end{eqnarray}

\z (see \cite{Wasow} below) where $m_i=1-\lfloor\beta_i\rfloor$,
($\lfloor\cdot\rfloor=$ integer part), $\bfC\in\CC^n$ is an arbitrary
vector of constants, and
$\tilde{\mathbf{y}}_\bfk=x^{-\bfk(\bfbet+\bfm)}
\sum_{l=0}^{\infty}\mathbf{a}_{\bfk;l} x^{-l}$ are formal power
series.  

\subsubsection{Summability} We give a brief  overview
of the results in \cite{Duke}.  The details of the proof follow the general
strategy presented in \S\ref{Ord1}. Let

\begin{eqnarray}
  \label{defW}
  \mathcal{W}=\left\{p\in\CC:p\ne k\lambda_i\,,\forall
k\in\NN,i=1,2,\ldots,n\right\}
\end{eqnarray}
(see Fig. 1) The directions $d_j=\{p:\arg(p)=\phi_j\}, j=1,2,\ldots,n$
are the {\em Stokes lines}.

We construct over $\mathcal{W}$ a surface $\mathcal{R}$,
consisting of homotopy classes of smooth curves in $\mathcal{W}$
starting at the origin, moving away from it, and crossing at most one
Stokes line, at most once (Fig. 1):

\begin{eqnarray}\label{defpaths}
{\mathcal{R}}:=\Big\{\gamma:(0,1)\mapsto \mathcal{W}:\ 
\gamma(0_+)=0;\ \frac{\mathrm{d}}{\mathrm{d}t}|\gamma(t)|>0;\ \arg(\gamma(t))\ \mbox{monotonic}\Big\}\cr
\end{eqnarray} 
Define $\mathcal{R}_1$ as the restriction of
$\mathcal{R}$ to $\arg(\gamma)\in(\psi_n-2\pi,\psi_2)$ where
$$\psi_n=\max\{-\pi/2,\phi_n-2\pi\}\ \text{and} \ \psi_2=\min\{\pi/2,\phi_2\}$$

\scalebox{0.5}[0.5]{\includegraphics{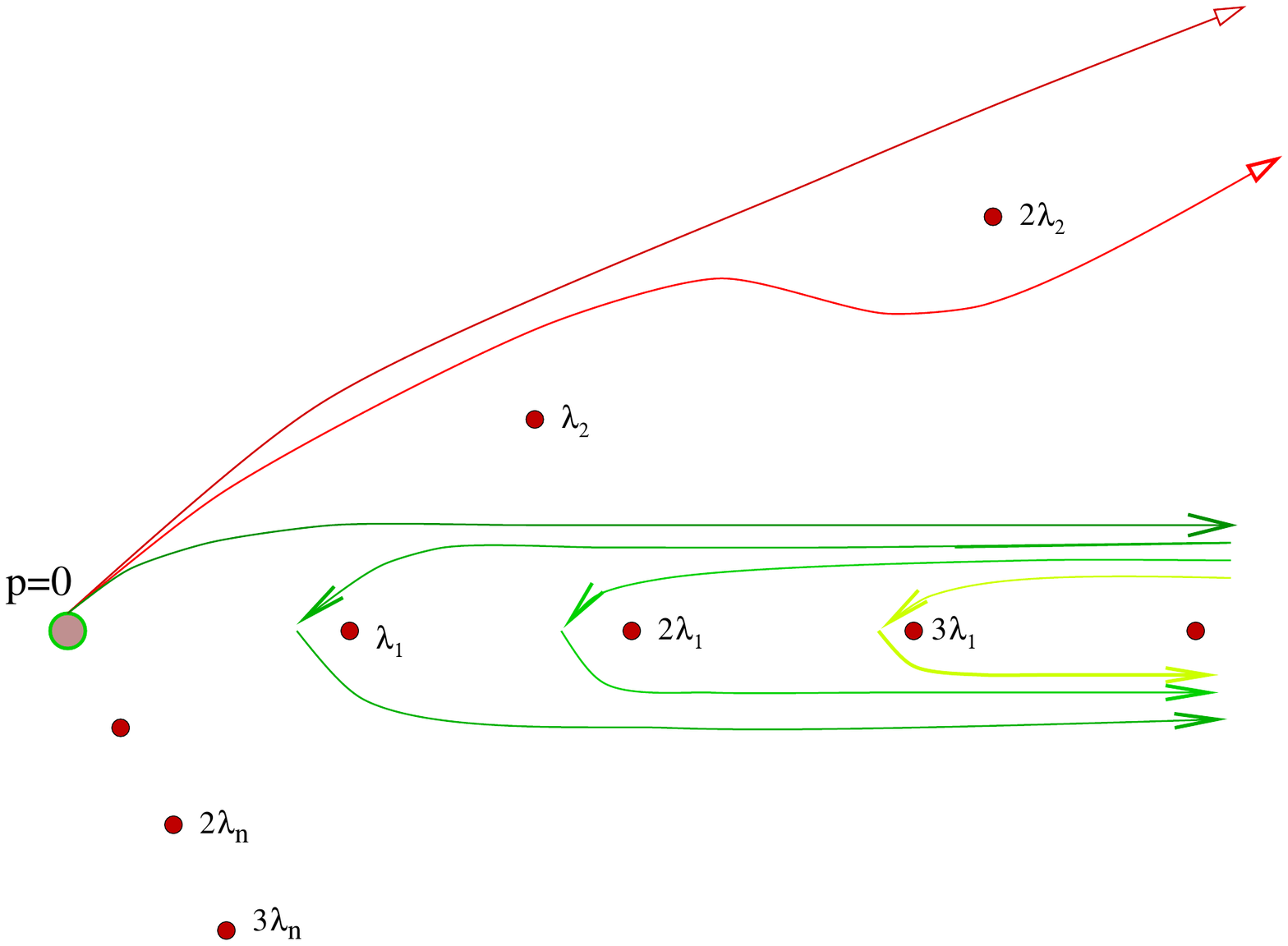}}

\centerline{{{Fig 1.} \emph{The paths
near $\lambda_2$ belong to $\mathcal{R}$.   }}}\nobreak
\centerline{\em The paths
near $\lambda_1$ relate to the balanced average}
\subsubsection{Singularities of $\bf Y_k$} The Borel transforms of $\bf \tilde{y}_k$ are analytic at zero and the only possible singularities 
are at multiples of the eigenvalues $\lambda_j$ of the linearized
system.  These singularities are ``regular'' in the sense that there
exist {\em convergent} local expansions at these singularities in
terms of powers and possibly logarithms. The following theorem makes
this statement precise.
 \begin{Theorem}[\cite{Duke}]\label{AS} (i) $\mathbf{Y}_0=\bor\tilde{\mathbf{y}}_0$
is analytic in $\mathcal{R}\cup \{0\}$. The singularities of $\mathbf{Y}_0$
(which are contained in the set
$\{l\lambda_j:l\in\NN^+,j=1,2,\ldots,n\}$) are described as follows.
For $l\in\NN^+$ and small $z$

 \begin{multline}
   \label{SY0}
   \mathbf{Y}_0^{\pm}(z+l\lambda_j)=\pm\Big[(\pm S_j)^l\ln(z)^{0,1}
   \mathbf{Y}_{l\mathbf{e}_j}(z)\Big]^{(lm_j)}+\mathbf{B}_{lj}(z)=\cr
   \Big[z^{l\beta_j'-1}\ln
   z^{0,1}\,\mathbf{A}_{lj}(z)\Big]^{(lm_j)}+\mathbf{B}_{lj}(z) \ (l=1,2,\ldots)
 \end{multline}

 \z where the power of $\ln(z)$ is one iff
 $l\beta_j\in\ZZ$, and $\mathbf{A}_{lj},\mathbf{B}_{lj}$ are analytic for small
 $z$. The functions $\mathbf{Y}_\bfk$ are, in addition, analytic at
 $p=l\lambda_j$, $l\in\NN^+$, iff, exceptionally,

\begin{eqnarray}\label{defSj}
S_j=r_j\Gamma(\beta'_j)\left(\mathbf{A}_{1,j}\right)_j(0)=0
\end{eqnarray}

\z where $r_j=1-\mathrm{e}^{2\pi \mathrm{i}(\beta'_j-1)}$
if $l\beta_j\notin\ZZ$ and $r_j=-2\pi \mathrm{i}$
otherwise. The $S_j$ are Stokes constants,
see theorem~\ref{Stokestr}.

(ii) $\mathbf{Y}_\bfk=\bor{\tilde{\mathbf{y}}}_\bfk$, $|\bfk|>1$, are analytic in
$\mathcal{R}\backslash
\{-\bfk'\cdot{\bflam}+\lambda_i:\bfk'\le\bfk,1\le i\le n\}$.  For
$l\in\NN$ and $p$ 
 near $l\lambda_j$, $j=1,2,\ldots,n$ there exist $\mathbf{A}=\mathbf{A}_{\bfk jl}$ and
$\mathbf{B}=\mathbf{B}_{\bfk jl}$ analytic at zero so that ($z$ is as above)

 \begin{multline}
   \label{SYK}
  \mathbf{Y}_\bfk^{\pm}(z+l\lambda_j)=
\pm\Big[(\pm S_j)^l{k_j+l\choose
     l}\ln(z)^{0,1}
   \mathbf{Y}_{\bfk+l\mathbf{e}_j}(z)\Big]^{(lm_j)}+l\mathbf{B}_{\bfk lj}(z)=\cr
  \Big[z^{\bfk\cdot{\bfbet}'+l\beta_j'-1}(\ln z)^{0,1}\,\mathbf{A}_{\bfk l
    j}(z)\Big]^{(lm_j)}+l\mathbf{B}_{\bfk l j}(z)\ (l=0,1,2,\ldots)
       \end{multline}

       \z where the power of $\ln z$ is $0$ iff $l=0$ or
       $\bfk\cdot{\bfbet}+l\beta_j-1\notin\ZZ$ and $\mathbf{A}_{\bfk 0
         j}=\mathbf{e}_j/\Gamma(\beta'_j)$.  Near
       $p\in\{-\bfk'\cdot{\bflam}:0\prec\bfk'\le\bfk\}$, (where
       $\mathbf{Y}_0$ is analytic) $\mathbf{Y}_\bfk,\,\bfk\ne 0$ have convergent
       Puiseux series.

\end{Theorem} 
\subsubsection{Averaging}\label{Ave} It is possible to take Laplace transforms
of $\bf Y_0$ along a ray avoiding $\RR^+$ from above or below. However
it can be checked that, assuming $\bf Y_0$ is singular at $p=1$ and
that the series $\tilde{y}_0$ has real coefficients, neither of them
would yield a real valued function.

 Let $\bor\tilde{\mathbf{y}}_\bfk$ be extended along $d_j$ by
the ``balanced average'' of analytic continuations

\begin{eqnarray}
\label{defmed}
\bor\tilde{\mathbf{y}}_\bfk=\mathbf{Y}_\bfk^{ba}=
\mathbf{Y}_\bfk^++\sum_{j=1}^{\infty}\frac{1}{2^j}\left(\mathbf{Y}_\bfk^{-}
-\mathbf{Y}_\bfk^{-({j-1})+}\right)
\end{eqnarray}

\z The sum above coincides with the one in which $+$ is exchanged with
$-$, accounting for the
reality-preserving property. Clearly, if $\mathbf{Y}_\bfk$ is analytic along
$d_j$, then the terms in the infinite sum vanish and
$\mathbf{Y}_\bfk^{ba}=\mathbf{Y}_\bfk$; we also let $\mathbf{Y}_\bfk^{ba}=\mathbf{Y}_\bfk$ if
$d\ne d_j$, where again $\mathbf{Y}_k$ is analytic. It follows from
(\ref{defmed}) and theorem~\ref{CEQ} below that the Laplace integral
of $\mathbf{Y}^{ba}_\bfk$ along $\RR^+$ can deformed into contours
as those depicted in Fig. 1, with weight $2^{-k}$ for a contour turning around
$(k+1)\lambda_1$. More generally, we consider the averages

\begin{eqnarray}
\label{defmedC0}
\bor_\alpha\tilde{\mathbf{y}}_\bfk=\mathbf{Y}_\bfk^{\alpha}=
\mathbf{Y}_\bfk^++\sum_{j=1}^{\infty}\alpha^j\left(\mathbf{Y}_\bfk^{-}
-\mathbf{Y}_\bfk^{-({j-1})+}\right)
\end{eqnarray}

\z and correspondingly

\begin{eqnarray}
  \label{defmedC}
(\mathcal{LB})_\alpha\tilde{\mathbf{y}}_\bfk:=\mathcal{L}\mathbf{Y}_\bfk^{\alpha}
\end{eqnarray}

\z With $\alpha\in\RR$, this represents the most general family of
averages of Borel summation formulas which commute with complex
conjugation, with the algebraic and analytic operations and have good
continuity properties \cite{Duke}.  The value $\alpha=1/2$ is special
in that it is the only one compatible with optimal truncation.

\begin{Note}
  For rank one ODEs the balanced average mentioned above can be shown
  to coincide with the ``median'' summation of Ecalle, who has
  introduced and studied a wide variety of averages \cite{EcalleMenous}, suitable
in more complicated settings.
\end{Note}

\begin{Theorem}\label{CEQ} (i) The branches of $(\mathbf{Y}_\bfk)_\gamma$ in $\mathcal{R}_1$
  have limits in a $C^*$-algebra of distributions,
  $\mathcal{D}'_{m,\nu}(\RR^+)\subset\mathcal{D}'$. Their Laplace transforms in
  $\mathcal{D}'_{m,\nu}(\RR^+)$ $\lap(\mathbf{Y}_\bfk)_\gamma$ exist
  simultaneously and with $x\in\mathcal{S}_x$ and for any $\delta>0$
  there is a constant $K$ and an $x_1$ large enough, so that for
  $\Re(x)>x_1$ we have $\left|\lap(\mathbf{Y}_\bfk\right)_\gamma(x)|\le
  K\delta^{|\bfk|}$.

In addition, $\mathbf{Y}_\bfk(p\mathrm{e}^{\mathrm{i}\phi})$ are continuous
in $\phi$ with respect to the $\mathcal{D}'_{m,\nu}$ topology,
(separately) on $[\psi_n-2\pi,0]$ and $[0,\psi_2]$.

If
$m>\max_i(m_i)$ and $l<\min_i |\lambda_i|$ then
$\mathbf{Y}_0(p\mathrm{e}^{\mathrm{i}\phi})$ is continuous in
$\phi\in[0,2\pi]\backslash\{\phi_i:i\le n\}$ in the
$\mathcal{D}'_{m,\nu}(\RR^+,l)$ topology and has (at most) jump
discontinuities for $\phi=\phi_i$. For each $\bfk$, $|\bfk|\ge 1$ and
any $K$ there is an $l>0$ and an $m$ such that
$\mathbf{Y}_k(p\mathrm{e}^{\mathrm{i}\phi})$ are continuous in
$\phi\in[0,2\pi]\backslash\{\phi_i; -\bfk'\cdot\bflam+\lambda_i:i\le n
,\bfk'\le\bfk\}$ in the $\mathcal{D}'_{m,\nu}((0,K),l)$ topology and
have (at most) jump discontinuities on the boundary.

(ii) The sum (\ref{defmed}) converges in $\mathcal{D}'_{m,\nu}$ (and
coincides with the analytic continuation of $\mathbf{Y}_\bfk$ when
$\mathbf{Y}_\bfk$ is analytic along $\RR^+$). For any
$\delta$ there is a large enough $x_1$ {\em
  independent of $\bfk$} so that  $\mathbf{Y}^{ba}_\bfk(p)$ with
$p\in\mathcal{R}_1$ are
Laplace transformable
for
$\Re(xp)>x_1$ and furthermore $|(\lap\mathbf{Y}^{ba}_\bfk)(x)|\le
\delta^{|\bfk|}$. In addition, if $d\ne\RR^+$, then for large $\nu$,
$\mathbf{Y}_\bfk\in L^1_\nu(d)$.

The functions
$\lap\mathbf{Y}_\bfk^{ba}$ are analytic for $\Re(xp)>x_1$. For any
$\bfC\in\CC^{n_1}$ there is an $x_1(\bfC)$ large enough so that the
sum

\begin{eqnarray}
  \label{soleqn}
  \mathbf{y}=\lap\mathbf{Y}_0^{ba}+\sum_{|\bfk|> 0}\bfC^{\bfk}\mathrm{e}^{-\bfk\cdot\bflam
    x}x^{-\bfk\cdot\bfbet}\lap\mathbf{Y}_\bfk^{ba}
\end{eqnarray}

\z converges uniformly for $\Re(xp)>x_1(\bfC)$, and $\mathbf{y}$ is a solution
of (\ref{eqor}).  When the direction
of $p$ is not the real axis then, by definition, $\mathbf{Y}^{ba}_\bfk=\mathbf{Y}_\bfk$,
$\mathcal{L}$ is the usual Laplace transform 
and (\ref{soleqn}) becomes

\begin{eqnarray}
  \label{soleqnpm}
  \mathbf{y}=\lap\mathbf{Y}_0+\sum_{|\bfk|> 0}\bfC^{\bfk}\mathrm{e}^{-\bfk\cdot\bflam
    x}x^{-\bfk\cdot\bfbet}\lap\mathbf{Y}_\bfk
\end{eqnarray}

In addition, $\lap\mathbf{Y}_\bfk^{ba}\sim \tilde{\mathbf{y}}_\bfk$ for large $x$
in the half plane $\Re(xp)>x_1$, for all $\bfk$, uniformly.

iii) The general
solution of (\ref{eqor}) that is asymptotic to $\tilde{\mathbf{y}}_0$ for
large $x$ along a ray
in $S_x$ can be equivalently written in  the form
(\ref{soleqn}) or as
\begin{eqnarray}
  \label{soleqnp}
  \mathbf{y}=\lap\mathbf{Y}_0^{\pm}+\sum_{|\bfk|> 0}\bfC^{\bfk}\mathrm{e}^{-\bfk\cdot\bflam
    x}x^{-\bfk\cdot\bfbet}\lap\mathbf{Y}_\bfk^{\pm}
\end{eqnarray}

\z for some $\bfC$ (depending on the solution and chosen form). With
the convention binding the directions of $x$ and $p$ and the
representation form being fixed the representation of a solution is
unique.

\end{Theorem}

\begin{Theorem}\label{RE}
i) For all $\bfk$ and $\Re(p)>j,\Im(p)>0$ as well
as in $\mathcal{D}'_{m,\nu}$ we have

\begin{eqnarray}
  \label{mainresur}
  \mathbf{Y}_{\bfk}^{\pm j\mp}(p)-\mathbf{Y}_{\bfk}^{\pm (j-1) \mp}(p) = (\pm S_1)^j\binom{k_1+j}{j}
  \left(\mathbf{Y}^\pm_{\bfk+j\mathbf{e}_1}(p-j)\right)^{(mj)}
\end{eqnarray}

\z and also, 

\begin{eqnarray}
  \label{thirdresu}
  \mathbf{Y}_{\bfk}^\pm=\mathbf{Y}_\bfk^{\mp}+\sum_{j\ge 1} {j+k\choose k}(\pm S_1)^{j}(\mathbf{Y}^\mp_{\bfk+j\mathbf{e}_1}(p-j))^{(mj)}
\end{eqnarray}

ii) {\em Local Stokes transition.}

\z Consider the expression of a fixed solution $\mathbf{y}$ of (\ref{eqor})
as a Borel summed transseries (\ref{soleqn}). As $\arg(x)$ varies,
(\ref{soleqn}) changes only through $\bfC$, and that change occurs
when the Stokes lines are crossed. We have, in the neighborhood
of $\RR^+$, with $S_1$
defined in (\ref{defSj}):

\begin{eqnarray}
  \label{microsto}
  \bfC(\xi)=\left\{\begin{array}{lll} \bfC^-=\bfC(-0)\qquad&\mbox{for
    $\xi<0$}\\ \bfC^0=\bfC(-0)+\frac{1}{2}S_1\mathbf{e}_1\qquad&\mbox{for
    $\xi=0$}\\ \bfC^+=\bfC(-0)+S_1\mathbf{e}_1\qquad&\mbox{for
    $\xi>0$}\end{array}\right.
\end{eqnarray}
\end{Theorem}

\begin{Remark}\label{R1} In view of (\ref{mainresur}) the different analytic
continuations of $\mathbf{Y}_0$ along paths crossing
$\RR^+$ at most once can be expressed in terms of
$\mathbf{Y}_{j\mathbf{e}_1}$. The most general formal solution of (\ref{eqor})
that can be formed in terms of $\mathbf{Y}_{j\mathbf{e}_j}$ with $j\ge 0$ is
(\ref{eqformgen,n}) with $C_1=\alpha$ arbitrary and $C_j=0$ for $j\ne
1$. Any true solution of (\ref{eqor}) based on such a transseries
is given in (\ref{soleqnp}) with 
$\bfC$ as above. Any average $\mathcal{A}\mathbf{Y}_0$  along paths
going forward in $\RR^+$ such that 
$\lap \mathcal{A}\mathbf{Y}_0$ is thus of the form (\ref{defmedC}).
\end{Remark}

\begin{Theorem}\label{Stokestr} Assume only $\lambda_1$ lies
  in the right half plane. Let $\gamma^{\pm}$ be two paths in the
  right half plane, near $\pm i\RR^+$ such
  that $|x^{-\beta_1+1}\mathrm{e}^{-x\lambda_1}|\rightarrow 1$ as
  $x\rightarrow\infty$ along $\gamma^{\pm}$. Consider the solution
  $\mathbf{y}$ of (\ref{eqor}) given in (\ref{soleqn}) with
  $\mathbf{C}=C\mathbf{e}_1$ and where the path of integration is
  $p\in\RR^+$. Then

\begin{eqnarray}\label{classicS}
\mathbf{y}=
(C\pm\frac{1}{2}S_1)\mathbf{e}_1 x^{-\beta_1+1}\mathrm{e}^{-x\lambda_1}(1+o(1))
\end{eqnarray}

\z for large $x$ along $\gamma^{\pm}$, where $S_1$ is the same as in
(\ref{defSj}), (\ref{microsto}).

\end{Theorem}

\begin{Proposition}\label{asymptrick} i)  Let $\mathbf{y}_1$ and $\mathbf{y}_2$ be solutions of (\ref{eqor})
so that $\mathbf{y}_{1,2}\sim\tilde{\mathbf{y}}_0$ for large $x$ in an
open sector $S$ (or in some direction
$d$);
then $\mathbf{y}_1-\mathbf{y}_2=\sum_{j}C_j\mathrm{e}^{-\lambda_{i_j}
  x}x^{-\beta_{i_j}}(\mathbf{e}_{i_j}+o(1))$ for some constants $C_j$, where the indices run over
  the eigenvalues $\lambda_{i_j}$ with the property $\Re(\lambda_{i_j}
  x)>0$ in $S$ (or $d$). 
If $\mathbf{y}_1-\mathbf{y}_2=o(\mathrm{e}^{-\lambda_{i_j}
    x}x^{-\beta_{i_j}})$ for all $j$, then $\mathbf{y}_1=\mathbf{y}_2$.
 
ii) Let $\mathbf{y}_1$ and $\mathbf{y}_2$ be solutions of
(\ref{eqor1}) and assume that $\mathbf{y}_1-\mathbf{y}_2$ has differentiable
asymptotics of the form 
$\mathbf{K}a\exp(-ax)x^b(1+o(1))$ with  $\Re(ax)>0$ and $\mathbf{K}\ne 0$, for large $x$.
Then $a=\lambda_i$ for some $i$.

iii) Let $\mathbf{U}_\bfk\in\mathcal{T}_{\{\cdot\}}$ for all $\bfk$,
$|\bfk|>1$.  Assume in addition that for
large $\nu$ there is a function $\delta(\nu)$ vanishing as
$\nu\rightarrow\infty$ such that

\begin{gather}
\label{cd1}
\sup_{\bfk}\delta^{-|\bfk|}\int_{d}\left|\mathbf{U}_{\bfk}(p)\mathrm{e}^{-\nu p}\right|\mathrm{d}|p|<K<\infty
\end{gather}

\z Then, if $\mathbf{y}_1,\mathbf{y}_2$ are solutions of (\ref{eqor}) in $S$ where in addition

\begin{eqnarray}
  \label{lapcond}
  \mathbf{y}_1-\mathbf{y}_2=\sum_{|\bfk|>1}
\mathrm{e}^{-\bflam\cdot\bfk x}x^{\bfm\cdot\bfk}\int_d\mathbf{U}_\bfk(p)
\exp(-xp)\mathrm{d}p
\end{eqnarray}

\z where $\bflam,x$ are as in (n6), then $\mathbf{y}_1=\mathbf{y}_2$, and
$\mathbf{U}_\bfk=0$ for all $\bfk$, $|\bfk|>1$.

\end{Proposition}

Given $\mathbf{y}$, the value of $C_i$ can change only when
$\xi+\arg(\lambda_i-\bfk\cdot\bflam)=0$, $k_i\in\NN\cup\{0\}$, i.e.
when crossing one of the (finitely many by (c1)) Stokes lines.  The
procedure is similar to the medianization \index{medianization}  proposed by \'Ecalle, but (due
to the structure of (\ref{eqor})) requires substantially fewer analytic
continuation paths. Resurgence \index{resurgence} relations are found and in addition we
provide a complete description, needed in applications, of the
singularity structure of the Borel transforms of $\tilde{\mathbf{y}}_\bfk$.

\subsubsection{More general equations amenable to first rank. Examples}\label{sec:ExaNorm} 

Many classical equations which are not presented in the form studied
in \S~\ref{Sysnonlin} can be brought to that form by elementary
transformations. We look at a few examples, to illustrate this and the
normalization process.

\z {\bf (1)}. The equation

\begin{equation}
  \label{erf}
  f'+2xf=1
\end{equation}

\z for large $x$ is not of the form (\ref{eqor1}). It can nevertheless
be brought to that form by simple changes of variables. The
normalizing change of variables of a given is most conveniently
obtained in the following way.  

{\em A transformation which brings the
  equation to its normal form also brings its transseries solutions to
  the form (\ref{transsf})}. 

 It is simpler to look for substitutions
with this latter property, and then the first step is to find the
transseries solutions of the equation.

At level zero, differentiation is contractive and thus, within power
series there is a unique solution of (\ref{erf}), obtained as a fizxed
point of
\begin{equation}
  \label{fixp}
  f=\frac{1}{x}-\frac{1}{x}f'
\end{equation}
Denoting this series by $f_1$ we look for further transseries solutions in the form $f=f_1+\delta$. We get

\begin{equation}
  \label{delta1}
    \delta'+2x\delta=0
\end{equation}
where we look for higher level transseries solutions (we choose to
ignore the fact that (\ref{delta1}) is solvable explicitly).  We
therefore take $\delta=e^w$ and look for level zero solutions.  We
have $w'+2x=0$ with the one parameter family of solutions
$w=-x^2/2+K$.  Therefore
\begin{equation}
  \label{tr2}
  f=f_1+Ce^{-x^2/2}
\end{equation}
The exponent is supposed to be linear in $x$ which suggests the change of variable $x^2=z$. In terms of $z$ we get
\begin{equation}
  \label{zz}
  g'+g=\frac{1}{2\sqrt{z}}
\end{equation}
which is not of the required form because of the noninteger power of $z$.
This can be easily take care of by the substitution $g=\sqrt{z} h$
which leads to
\begin{equation}
  \label{fin3}
  h'+\Big(1+\frac{1}{2z}\Big)h=\frac{1}{2z}
\end{equation}

{\bf (2)} A simple Schrodinger equation in one dimension (the harmonic
oscillator).

\begin{equation}
  \label{eq:Schr}
y''-x^2y=\lambda y  
\end{equation}

A convenient way to find transseries solutions is to use a WKB-like
procedure. Substituting $y(x)=e^{g(x)}$ in (\ref{eq:Schr}) one obtains

$$g'=\pm \sqrt{\lambda+x^2-g''}$$

\z which is a contractive equation in the space
$\tilde{\mathcal{T}}^{[0]}$, cf. \S~\ref{deg1}. It follows 
in particular that $\tilde{g}'=\pm x+O(1/x)$ and thus $\ln
y(x)=\pm\frac{1}{2}x^2+O(\ln x)$. Since for (\ref{eqor1})
the exponentials have linear exponents, the natural variable
of (\ref{eq:Schr}), for our purpose, is $t=\frac{1}{2} x^2$. In this
variable, (\ref{eq:Schr}) becomes

\begin{equation}
  \label{eq:Schr,norm}
h''+\frac{1}{2t}h'-\left(1+\frac{\lambda}{2t}\right)h=0
\end{equation}

\z which in vector form, with $y_1=h+h', y_2=h-h'$, reads

\begin{equation}
  \label{eq:SchrSyst}
  \mathbf{y}'=
\left\{\begin{pmatrix}
  1&0\cr 0&-1
\end{pmatrix}+\frac{1}{4t}
\begin{pmatrix}
  \lambda-1&\lambda+1\cr
1-\lambda&-1-\lambda 
\end{pmatrix}\right\}
\mathbf{y}=\left(\hat{\Lambda}+\frac{1}{t}\hat{C}\right)
\mathbf{y}
\end{equation}

The matrix $\hat{C}$ can (always) be diagonalized, by taking
$\mathbf{y}=(I+t^{-1}\hat{S})\mathbf{Y}$ for a suitable matrix
$\hat{S}$, and the diagonalized matrix is
$\hat{B}=\mbox{diag}(C_{11},C_{22})$.  Indeed, the equation for
$\hat{S}$ is
$\hat{\Lambda}\hat{S}-\hat{S}\hat{\Lambda}=\hat{B}-\hat{C}$, which can
be solved whenever the eigenvalues of $\hat{\Lambda}$ are distinct.  The
normal form of (\ref {eq:SchrSyst}) is thus

\begin{equation}
  \label{eq:SchrSystN}
  \mathbf{y}'=\left(\hat{\Lambda}+\frac{1}{t}\hat{C}\right)
\mathbf{y}+t^{-2}\hat{g}(t^{-1})\mathbf{y}
\end{equation}

\z with $\hat{g}$ analytic. 

{\bf (3)} A nonintegrable case of Abel's equation

\begin{gather}
  \label{(*)}
 w'=w^3-z 
\end{gather}
{\em Power series solutions}. As before, (\ref{(*)}) written in the form

$$w=(w'+z)^{1/3}$$

\z is contractive in $\tilde{\mathcal{T}}^{[0]}$. By iteration,  one
obtains the power series formal solutions
$\tilde{w}_0=Az^{1/3}(1+\sum_{k=1}^{\infty}{\tilde{w}}_{0,k}z^{-5k/3})$
($A^3=1$).

\z {\em General transseries solutions of (\ref{(*)})}. Once we determined
$\tilde{w}_0$ we look for possible small corrections to $\tilde{w}_0$;
since these are usually exponentially small, formal WKB 
is again useful. We substitute $w=\tilde{w}_0+e^g$. The 
equation for $g$

$$\tilde{g}=C+\mathcal{P}\left(3\tilde{w}_0^2+3\tilde{w}_0 e^{\tilde{g}}
  +e^{2\tilde{g}}\right)$$

\z is contractive in any sector where $\Re(\tilde{w}_0)<0$ and in this
case $e^{\tilde{g}}\propto z^{2/3}\exp\left(\frac{9}{5}A^2z^{5/3}\right)$

Since the exponentials in a transseries solution of a normalized system
have {\em linear} exponent, with negative real part, the independent
variable should be $x=-(9/5)A^2z^{5/3}$ and $\Re(x)>0$. Then
$\tilde{w}_0=x^{1/5}$ $\sum_{k=0}^{\infty} w_{0;k} x^{-k}$, which
compared to (\ref{transsf}) suggests the change of dependent variable
$w(z)=Kx^{1/5}h(x)$. Choosing for convenience $K=A^{3/5}(-135)^{1/5}$
yields

\begin{gather}
  \label{tr12}
h'+\frac{1}{5x}h+3h^3-\frac{1}{9}=0
\end{gather}

\z The next step is to achieve leading behavior $O(x^{-2})$. This is
easily done by subtracting out the leading behavior of $h$ (which can be
found by maximal balance, as above). With $h=y+1/3-x^{-1}/15$ we
get the normal form

\begin{align}
  \label{eqfu}
y'=-\left(1-\frac{1}{5x}\right)y+g(y,x^{-1})
\end{align}

\z where
\begin{equation}
  \label{eq:defg}
  g(y,x^{-1})=-3(y^2+y^3)+\frac{3y^2}{5x}-\frac{1}{15x^2}-\frac{y}{25x^2}+\frac{1}{3^2
  5^3 x^3}
\end{equation}

{\bf (3)} The Painlev\'e equation P1.

\begin{gather}
  \label{eP1}
\frac{d^2y}{dz^2}=6y^2+z
\end{gather}

The transformations needed to normalize
(\ref{eP1}) are derived in the same way as in Example 1.  After the
change of variables
\begin{gather*}
x=\frac{(-24z)^{5/4}}{30};\ y(z)=\sqrt{\frac{-z}{6}}\left(1-\frac{4}{25x^2}+h(x)\right)
\end{gather*}
P1 becomes 

\begin{gather}\label{eqp1n}
h'' +\frac{1}{x}h'-h-\frac{1}{2}h^2-\frac{392}{625x^4}=0
\end{gather} 

\z Written as a system, with $\mathbf{y}=(h,h')$ this equation satisfies
the assumptions in \S~\ref{Sysnonlin} with $\lambda_{1,2}=\pm 1$,
$\alpha_{1,2} =-1/2$, and then $\xi(x)=C\erm^{-x}x^{-1/2}$.

{\bf (4)} The Painlev\'e equation P2. 

This equation reads:

\begin{eqnarray}
  \label{eqP2}
  y''=2y^3+xy+\alpha
\end{eqnarray}

\z This example also shows that for a given equation distinct solution
manifolds associated to distinct asymptotic behaviors may lead to
different normalizations. After the change of variables

$$x=(3t/2)^{2/3};\ \ \
y(x)=x^{-1}(t\,h(t)-\alpha)$$

\z one obtains the normal form equation
\begin{align}
  \label{eq:p2n0}
  h''+\frac{h'}{t}-\left(1+\frac{24\alpha^2+1}{9t^2}\right)h-\frac{8}{9}h^3+
\frac{8\alpha}{3t}h^2+\frac{8(\alpha^3-\alpha)}{9t^3}=0
\end{align}

Distinct  normalizations (and sets of solutions)
are provided by 

$$x=(At)^{2/3};\ \
y(x)=(At)^{1/3}\left(w(t)-B+\frac{\alpha}{2At}\right)$$

\z if $A^2=-9/8,B^2=-1/2$. In this case, 

\begin{multline}
  w''+\frac{w'}{t}+w\left(1+\frac{3B\alpha}{tA}-
\frac{1-6\alpha^2}{9t^2}\right)w\\
-\left(3B-\frac{3\alpha}{2tA}\right)w^2+w^3+
\frac{1}{9t^2}\left(B(1+6\alpha^2)-t^{-1}\alpha(\alpha^2-4)
\right)
\end{multline}

\z so that

$$\lambda_1=1, \alpha_1=-\frac{1}{2}-\frac{3}{2}\frac{B\alpha}{A}$$

The first normalization applies for the manifold of solutions such
that $y\sim-\frac{\alpha}{x}$ (for $\alpha=0$ $y$ is exponentially
small and behaves like an Airy function) while the second one
corresponds to $y\sim -B-\frac{\alpha}{2}x^{-3/2}$.
\subsection{Difference equations and PDEs}
\subsubsection{Difference equations}  Transseries and
Borel summation techniques can be successfully used for difference equations.
Consider difference systems of equations which can be brought to
the form
\begin{equation}
  \label{eq:Bra1}
  \mathbf{x}(n+1)=\hat\Lambda\left(I+ \frac{1}{n}\hat{A}\right)\mathbf{x}(n) +\mathbf{g}(n,\mathbf{x}(n))
\end{equation}
where $\hat\Lambda $ and $\hat{A} $ are constant coefficient matrices, $\mathbf{g}$ is
convergently given for small $\mathbf{x}$ by
\begin{equation}
  \label{eq:defg1}
\mathbf{g}(n,\mathbf{x})=\sum_{\mathbf{k}\in\NN^m}
\mathbf{g}_\mathbf{k}(n)\mathbf{x}^\mathbf{k} 
\end{equation}
with $\mathbf{g}_\mathbf{k}(n)$ analytic in $n$ at infinity and
\begin{equation}
  \label{eq:condg}
  \mathbf{g}_\mathbf{k}(n)=O(n^{-2}) \ \mbox{as }n\rightarrow\infty,\
  \mbox{if }\sum_{j=1}^m k_j\le 1
\end{equation}
 under nonresonance conditions: Let
$\boldsymbol{\mu}=(\mu_1,...,\mu_n)$ and $\mathbf a=(a_1,...,a_n)$ where
$e^{-\mu_k}$ are the eigenvalues of $\hat{\Lambda}$ and the $a_k$ are
the eigenvalues of $\hat{A}$.  Then the nonresonance condition is
\begin{equation}
  \label{eq:nonres}
  \left(\mathbf{k}\cdot \boldsymbol\mu=0 \mod 2\pi i \ \ \mbox{with}\ \ \mathbf{k} \in\ZZ^{m_1}
  \right)
\Leftrightarrow \mathbf{k}=0.
\end{equation}
We consider the solutions of (\ref{eq:Bra1}) which are small as $n$
becomes large.  \ Braaksma and Kuik \cite{Braaksma} \cite{Kuik}
showed that the recurrences (\ref{eq:Bra1}) posess $l$-parameter
transseries solutions of the form
\begin{equation}
  \label{eq:transszf}
  \tilde{{\bf x}}(t)
:=\sum_{{\bf k}\in\NN^{m}}{\bf C}^{\bf k} e^{-{\bf k}\cdot
  \boldsymbol{\mu} t}
t^{{\bf k}\cdot {\bf a}} \tilde{{\bf x}}_\bfk(t)
\end{equation}
with $t=n$ where $ \tilde{{\bf x}}_\bfk(n)$ are formal power series in
powers of $n^{-1}$ and $l\le m$ is chosen such that, after reordering
the indices, we have $\Re(\mu_j)>0$ for $1\le j\le l$.
 
 It is shown in \cite{Braaksma}, \cite{Kuik} that these transseries are generalized
 Borel summable in any direction and Borel summable in all except $m$
 of them and that
\begin{equation}
  \label{eq:transsSn}
  {\bf x}(n)
=\sum_{{\bf k}\in\NN^{l}}{\bf C}^{\bf k} e^{-{\bf k}\cdot
  \boldsymbol{\mu} n}
n^{{\bf k}\cdot {\bf a}}{\bf x}_\bfk(n)
\end{equation} 
is a solution of (\ref{eq:Bra1}), if $n>y_0$, $t_0$ large enough.
\subsubsection{PDEs}

\subsubsection{Existence, uniqueness, regularity, asymptotic behavior}
For partial differential equations with analytic coefficients which
can be transformed to equations in which the differentiation order in
a distinguished variable, say time, is no less than the one with
respect to the other variable(s), under some other natural
assumptions, Cauchy-Kowalevski theory (C-K) applies and gives
existence and uniqueness of the initial value problem. A number of
evolution equations do not satisfy these assumptions and even if
formal power series solutions exist their radius of convergence is
zero.  
The paper \cite{Invent2} provides a C-K type theory in such
cases, providing existence, uniqueness and regularity of the
solutions.  Roughly, convergence is replaced by Borel summability,
although the theory is more general.

Unlike in C-K, solutions of nonlinear evolution equations develop
singularities which can be more readily studied from the local
behavior near $t=0$, and this is useful in determining and proving
spontaneous blow-up.

In the following, $\partial_{\bf x}^{\bf j} \equiv
\partial_{x_1}^{j_1} \partial_{x_2}^{j_2} ...\partial_{x_d}^{j_d}$,
$|{\bf j}|= j_1 + j_2 + ... + j_d$, ${\bf x}$ is in a poly-sector
$\mathcal S=\{\mathbf x:|\arg x_i|<\frac{\pi}{2}+\phi;|\mathbf x|>a\}$
in $\CC^d$ where $\phi<\frac{\pi}{2n}$, ${\bf g} \left ( {\bf x}, t,
  \{\mathbf y_{\bf j}\}_{|{\bf j}|=0}^{n-1} \right )$ is a function
analytic in $\{\mathbf y_{\bf j}\}_{|{\bf j}|=0}^{n-1}$ near $\bf 0$
vanishing as $|\bf x|\to \infty$.  The results in \cite{Invent2} hold
for $n$-th order nonlinear {\em quasilinear} partial differential
equations of the form
\begin{equation}
  \label{geneq}
{\bf u}_t + \mathcal{P}(\partial_{\bf x}^{\bf j}){\bf u}+{\bf g}
  \left ( {\bf x}, t, \{\partial_{\bf x}^{{\bf j}} {\bf u}\} \right )
  =0
\end{equation}
\z where $\mathbf u\in\CC^m$, for large $|\bf x|$ in $S$.
Generically, the constant coefficient operator $\mathcal{P}(\partial_{\bf
  x}) $ in the linearization of $\mathbf g(\infty,t,\cdot)$ is
diagonalizable. It is then taken to be diagonal, with eigenvalues
$\mathcal{P}_j$.  $\mathcal{P}$ is subject to the requirement that for
all $j\le m$ and ${\bf p}\ne 0$ in $\mathbb{C}^d$ with $|\arg p_i|\le
\phi $ we have
\begin{equation}
  \label{conecond}
  \Re\mathcal{P}^{[n]}_j (-{\bf p} )> 0
\end{equation}
where $\mathcal{P}^{[n]}(\partial_{\bf x})$ is the principal symbol of
$\mathcal{P}(\partial_{\bf x})$. Then the following holds. (The
precise conditions and results are given in \cite{Invent2}.)
\begin{Theorem}[large $|\bf x|$ existence] 
 \label{T1} Under the assumptions above, for any $T>0$   (\ref{geneq}) 
 has a unique solution ${\bf u}$ that for $t\in [0,T]$ is $O(|\mathbf x|^{-1})$  and analytic in $\mathcal S$.
\end{Theorem}

Determining asymptotic properties of solutions of PDEs is
substantially more difficult than the corresponding question for ODEs.
Borel-Laplace techniques however provide a very efficient way to
overcome this difficulty.  The paper shows that formal series
solutions are actually Borel summable, a fortiori asymptotic, to
actual solutions.
\begin{Condition}\label{Cond 2}
  The functions ${\bf b}_{{\bf q}, {\bf k}} ({\bf x}, t)$ and ${\bf r}
  ({\bf x}, t)$ are analytic in
  $(x_1^{-\frac{1}{N_1}},...,x_d^{-\frac{1}{N_d}})$ for large $\bf|x|$
  and some $N\in\NN$.
\end{Condition}

\begin{Theorem}\label{TrB}
  If Condition~\ref{Cond 2} and the assumptions of Theorem~\ref{T1}
  are satisfied, then the unique solution $\mb f$ found there
  can be written as
  \begin{equation}
    \label{acc3}
    \mathbf{f}(\mathbf{x},t)=
\int_{{\RR^+}^d}e^{-\mathbf{p}\cdot\mathbf{x}^{\frac{n}{n-1}} }{\mathbf{F}_1}(\mathbf{p},t)d\mathbf{p}
  \end{equation}
  where ${\mathbf{F}_1}$ is (a)  analytic at zero in
  $(p_1^{\frac{1}{nN_1}},...,p_d^{\frac{1}{nN_d}})$; (b) analytic in $\bf p\ne 0$ in the
  poly-sector $|\arg p_i|<\frac{n}{n-1}\phi+\frac{\pi}{2(n-1)}$, $i\le
  d$; and (c) exponentially bounded in the latter poly-sector.
\end{Theorem}
Existence and asymptoticity of the formal power series follow as a
corollary, using Watson's Lemma. 

Eariler, Borel summability has been shown for heat equation by Lutz,
Miyake and Sch\"afke \cite{Lutz} and generalized to linear PDEs with
constant coefficients by Balser \cite{Balser} and in special classes
of higher order nonlinear PDEs in \cite{CPAM}.
\subsection{More general irregular singularities and multisummability} The
simple normalizing procedure described above does not always work;
although this situation is in some sense nongeneric it plays an
important role in some applications. In some instances mixed type
transseries; for instance, by taking $e^{-x}Ei(x)$ which is a solution
of a rank one differential equation as the rhs of (\ref{erf}) we get

\begin{equation}
  \label{erfmix}
  f'+2xf=e^{-x}Ei(x)
\end{equation}
which can be brought to a second order meromorphic equation by using
the equation of $e^{-x}Ei(x)$, i.e. by adding (\ref{erfmix}) to its
derivative,
\begin{equation}
  \label{erfmix2}
  Df:=f''+(2x+1)f'+2(1+x)f=\frac{1}{x}
\end{equation}
The power series solution can be obtained contractively as in the
previous examples. For the complete transseries solution, the
substitution $f=e^w$ now yields
$$\frac{1}{2}{w'}^2+xw'+x=0$$
which has two solutions of level zero,
\begin{eqnarray}
  \label{degsolns}
  & \displaystyle w_1=-x^2+x+\frac{1}{2}\ln x+K_1-\frac{1}{2x}-...\nonumber \\
  & \displaystyle w_2=-x-\frac{1}{2}\ln x+K_2+\frac{1}{2x}+...
\end{eqnarray}
and thus 
\begin{equation}
  \label{gentr}
  f=f_1+C_2 e^{-x^2+x} x^{\frac{1}{2}}f_2+ C_3 e^{-x}x^{-\frac{1}{2}}f_3
\end{equation}
with $f_1, f_2, f_3$ power series. 

Clearly, no change of independent variable brings (\ref{gentr}) to the
form (\ref{transsf}). On the other hand, neither of variable $x$ nor
$x^2-x$ (equivalently, $x^2$) can be used for Borel summation:

(i) Borel transform in the variable $x^2$ (or $x^2-x$) yields a
function that has still a {\em divergent} power series at the origin.
Indeed the the coefficients $c_n$  of the power series
$f_1$, which satisfy the recurrence relation
\begin{equation}
  \label{recrel}
  c_{n+1}=(n-1)\Big[c_{n}+\frac{1}{2}c_{n-1}-\frac{1}{2}(n-2)c_{n-2}\Big]
\end{equation}
(with $c_0=0=c_1=0, c_2=\frac{1}{2}$) grow like $n!$, while a Borel
transform in $x^2$ only decreases the growth of the coefficients by a
factor of, roughly, $(n!)^{-1/2}$. At a closer look it can be shown
seen that the factorial growth of the coefficients and the presence of
the term $e^{-x}$ are interrelated.

(ii) It can be checked that in the variable $x$, the Borel transform
$F$ of the power series solution $f$ of (\ref{erfmix2}) is not Laplace
transformable because of faster than exponential growth.  Indeed, $F$
satisfies the equation
\begin{equation}
  \label{Borelt1}
 2F'-pF =\frac{1}{1-p};\ \ \ F(0)=0
\end{equation}
The superexponential growth, in turn, can be related to the presence
of $e^{-x^2+x}$ in the transseries. 

\subsubsection{}\label{751} This suggests decomposing the formal power series solution into two
parts, one summable in the variable $x$, the other one in the variable
$x^2$. We can adjust the parameters $a$ and $b$ so that

(i) The initial condition $c_0=0,c_1=-a/2,
c_3=-\frac{a}{4}+\frac{b}{2}$ in the recurrence (\ref{recrel}) of the
equation $Df=-a-bx^{-2}$ implies at most  $(n!)^{-1/2}$ growth of $c_n$
and at the same time,

(ii) The Borel transform
of equation $Df=a+x^{-1}+bx^{-2}$ has a solution with subexponential
growth. 

The formal series solution of $Df=x^{-1}$ is thus decomposed in a sum
of separately Borel summable series in the variables $x^2$ and $x$
respectively.  That this (rather ad-hoc) procedure works is a
reflection of general theorems in multisummability.

\subsubsection{Multisummability} A powerful and general technique, that of {\em acceleration and
  multisummability} introduced by \'Ecalle, see \cite{Eca84} and
\cite{EcalleNATO}, adequately deals with mixed divergences in very
wide generality.  The procedure is universal and works in the same way
for very general systems.

In the case of solutions of meromorphic differential equations, mixed
types of diveregence relate to the presence of exponential terms with
exponents of different powers. Then multisummation consists in Borel
transform with respect to the lowest power of $x$ in the exponents of
the transseries, usual summation $\mathcal{S}$, a sequence of
transformations called accelerations (which mirror in Borel space the
passage from one power in the exponent to the immediately larger one)
followed by a final Laplace transform in the largest power of $x$.
More precisely (\cite{EcalleNATO}):
\begin{equation}
  \label{multisum}
  \mathcal{L}_{k_1}\circ \mathcal{A}_{k_2/k_1}\circ \cdots\circ \mathcal{A}_{k_{q}/k_{q-1}} \mathcal{S}\mathcal{B}_{k_q}
\end{equation}
where $ (\mathcal{L}_{k}f)(x)=(\mathcal{L}f)(x^k)$, $\mathcal{B}_{k}$ is the
formal inverse of $\mathcal{L}_{k}$, $\alpha_i\in (0,1)$ and the
acceleration operator $\mathcal{A}_{\alpha}$ is formal the image, in
Borel space, of the change of variable from $x^{\alpha}$ to $x$ and is
defined as
\begin{equation}
  \label{accel}
  \mathcal{A}_{\alpha}\phi=\int_0^{\infty}C_{\alpha}(\cdot,s)\phi(s)ds
\end{equation}
and where, for $\alpha\in(0,1)$,  the kernel $C_{\alpha}$ is defined as
\begin{equation}
  \label{accel2}
  C_{\alpha}(\zeta_1,\zeta_2):=\frac{1}{2\pi i}\int_{c-i\infty}^{c+i\infty}e^{\zeta_2 z-\zeta_1 z^{\alpha}}dz
\end{equation}
where we adapted the notations in \cite{Balser} to the fact that the
formal variable is large.  In our example, $q=2, k_2=1, k_1=2$.

In \cite{[Br]} Braaksma proved of multisummability of series solutions
of general nonlinear meromorphic ODEs using Borel transforms in the
spirit of Ecalle'e theory.

\begin{Note}
  (i) Multisummability of type (\ref{multisum}) can be equivalently
  characterized by decomposition of the series into terms which are
  ordinarily summable after changes of independent variable of the
  form $x\to x^\alpha$.  This is shown in \cite{Balser} where it is
  used to give an alternative proof of multisummability of series
  solutions of meromorphic ODEs, closer to the cohomological point of
  view of Ramis and Sibuya.
  
  (ii) More general multisummability is described by Ecalle
  \cite{EcalleNATO}, allowing, among others, for stronger than
  power-like acceleration. This is relevant to more general
  transseries equations.
  
  (iii) We expect that the finite generation property of transseries
  allows for the cohomological approach (i) combined with (ii) to
  prove multisummability of transseries solutions of more general
  contractive equations, in the sense of Theorem~\ref{PFp}. This would
  give a rigorous proof of closure of analyzable functions under 
  all natural operations in the sense of Ecalle.
\end{Note}

\end{document}